\documentclass[a4paper,11pt]{article}

\usepackage{comment} 
\usepackage[latin1]{inputenc}
\usepackage[english]{babel}   
\usepackage{amssymb}
\usepackage{amsmath}
\usepackage{amsthm}
\usepackage{multicol}
\usepackage{graphicx}
\usepackage{color}
\usepackage[T1]{fontenc}
\usepackage{yfonts}
\usepackage{placeins}
\usepackage{array}
\usepackage{dsfont}
\usepackage{stmaryrd}
\usepackage{verbatim}
\usepackage{caption}

\theoremstyle{plain}
\newtheorem{theorem}{Theorem}
\newcommand{\bthm}{\begin{theorem}}
\newcommand{\ethm}{\end{theorem}}

\newtheorem{thmi}{Theorem}
\newcommand{\bthmi}{\begin{thmi}}
\newcommand{\ethmi}{\end{thmi}}

\newtheorem{cori}[thmi]{Corollary}
\newcommand{\bcori}{\begin{cori}}
\newcommand{\ecori}{\end{cori}}

\newtheorem{mthm}{Theorem}
\newcommand{\bmthm}{\begin{mthm}}
\newcommand{\emthm}{\end{mthm}}
\newtheorem{mcor}[mthm]{Corollary}
\newcommand{\bmcor}{\begin{mcor}}
\newcommand{\emcor}{\end{mcor}}

\newtheorem*{conj}{Conjecture}
\newcommand{\bconj}{\begin{conj}}
\newcommand{\econj}{\end{conj}}

\newtheorem*{question}{Question}
\newcommand{\bq}{\begin{question}}
\newcommand{\eq}{\end{question}}

\newtheorem*{thn}{Theorem}
\newcommand{\bthn}{\begin{thn}}
\newcommand{\ethn}{\end{thn}}

\newtheorem{exo}{Exercise}
\newcommand{\bex}{\begin{exo}}
\newcommand{\eex}{\end{exo}}

\newtheorem{sol}{Solution}
\newcommand{\bsol}{\begin{sol}}
\newcommand{\esol}{\end{sol}}

\newtheorem{pro}[theorem]{Proposition}
\newcommand{\bpro}{\begin{pro}}
\newcommand{\epro}{\end{pro}}

\newtheorem{cor}[theorem]{Corollary}
\newcommand{\bcor}{\begin{cor}}
\newcommand{\ecor}{\end{cor}}

\newtheorem{lem}[theorem]{Lemma}
\newcommand{\blem}{\begin{lem}}
\newcommand{\elem}{\end{lem}}

\theoremstyle{definition}

\newtheorem{defi}{Definition}
\newcommand{\bdf}{\begin{defi}}
\newcommand{\edf}{\end{defi}}

\newtheorem*{rmq}{Remark}
\newcommand{\brq}{\begin{rmq} \upshape}
\newcommand{\erq}{\end{rmq}}

\newtheorem*{exe}{Example}
\newcommand{\bexe}{\begin{exe} \upshape}
\newcommand{\eexe}{\end{exe}}

\newcommand{\bp}{\begin{proof}}
\newcommand{\ep}{\end{proof}}

\newcommand{\beq}{\begin{eqnarray*}}
\newcommand{\eeq}{\end{eqnarray*}}

\newcommand{\beqn}{\begin{equation}}
\newcommand{\eeqn}{\end{equation}}

\newcommand{\ben}{\begin{enumerate}}
\newcommand{\een}{\end{enumerate}}
\newcommand{\bit}{\begin{itemize} \renewcommand{\labelitemi}{$\bullet$} }
\newcommand{\eit}{\end{itemize}}

\newenvironment{Figure}
  {\par\medskip\noindent\minipage{\linewidth}}
  {\endminipage\par\medskip}

\newcommand{\bfg}{
\begin{figure}[H]
\begin{center}}
\newcommand{\efg}{
\end{center}
\end{figure}
\FloatBarrier}

\newcommand{\df}{\emph}

\newcommand{\RRR}{\mathbb{R}}
\newcommand{\QQQ}{\mathbb{Q}}
\newcommand{\NNN}{\mathbb{N}}
\newcommand{\ZZZ}{\mathbb{Z}}
\newcommand{\CCC}{\mathbb{C}}

\newcommand{\PPP}{\mathbb{P}}

\newcommand{\KKK}{\mathbb{K}}

\newcommand{\HHH}{\mathbb{H}}

\newcommand{\FFF}{\mathbb{F}}

\newcommand{\bs}{\symbol{92}}

\newcommand{\ov}{\overline}

\renewcommand{\dim}{\operatorname{dim}}
\newcommand{\Card}{\operatorname{Card}}
\newcommand{\Isom}{\operatorname{Isom}}

\newcommand{\CAT}{\operatorname{CAT}}
\newcommand{\Sp}{\operatorname{Sp}}

\newcommand{\bord}{\partial_\infty}
\newcommand{\Mor}{\operatorname{Mor}}

\newcommand{\sa}{\sphericalangle}

\newcommand{\st}{\operatorname{st}}
\newcommand{\Gr}{\mathcal{G}r}

\newcommand{\PU}{\operatorname{PU}}
\newcommand{\PSp}{\operatorname{PSp}}

\newcommand{\lb}{\llbracket}
\newcommand{\rb}{\rrbracket}
\newcommand{\ra}{\rightarrow}
\newcommand{\lra}{\longrightarrow}
\newcommand{\lmapsto}{\longmapsto}
\newcommand{\lhookra}{\lhook\joinrel\longrightarrow}
\newcommand{\ral}[1]{\underset{#1}{\longrightarrow}}

\newcommand{\liml}{\lim\limits}

\renewcommand{\geq}{\geqslant}
\renewcommand{\leq}{\leqslant}

\newcommand{\id}{\operatorname{id}}

\newcommand{\SL}{\operatorname{SL}}
\newcommand{\SO}{\operatorname{SO}}
\newcommand{\SU}{\operatorname{SU}}
\newcommand{\PGL}{\operatorname{PGL}}
\newcommand{\PSL}{\operatorname{PSL}}

\newcommand{\PO}{\operatorname{PO}}

\newcommand{\Flats}{\operatorname{Flats}}
\newcommand{\Cham}{\operatorname{Cham}}
\newcommand{\Ap}{\operatorname{Ap}}
\newcommand{\<}{\langle}
\renewcommand{\>}{\rangle}

\makeatletter
\def\Ddots{\mathinner{\mkern1mu\raise\p@
\vbox{\kern7\p@\hbox{.}}\mkern2mu
\raise4\p@\hbox{.}\mkern2mu\raise7\p@\hbox{.}\mkern1mu}}
\makeatother

\makeatletter
\def\maketitles{%
  \null
  \thispagestyle{empty}%
  \vfill
  \begin{center}\leavevmode
    \normalfont
    {\LARGE \@title\par}%
    \vskip 1.2cm
    {\large \@author\par}%
    \vskip 1.2cm
    {\large \@subtitle\par}%
    \vskip 0.8cm
    {\large \@date\par}%
  \end{center}%
  \vfill
  \null
  \cleardoublepage
  }

\def\date#1{\def\@date{#1}}
\def\author#1{\def\@author{#1}}
\def\title#1{\def\@title{#1}}
\def\subtitle#1{\def\@subtitle{#1}}
\makeatother

\title{Visual limits of maximal flats\\in symmetric spaces and Euclidean buildings}
\author{Thomas Haettel}
\date{\today}

\begin{document}

\selectlanguage{english}

\maketitle

\begin{abstract}
Let $X$ be a Riemannian symmetric space of non-compact type or a locally finite, strongly transitive Euclidean building, and let $\bord X$ denote the geodesic boundary of $X$. We reduce the study of visual limits of maximal flats in $X$ to the study of limits of apartments in the spherical building $\bord X$: this defines a natural, geometric compactification of the space of maximal flats of $X$. We then completely determine the possible degenerations of apartments when $X$ is of rank $1$ or associated to a classical group of rank $2$ or to $\PGL(4)$. In particular, we exhibit remarkable behaviours of visual limits of maximal flats in various symmetric spaces of small rank and the surprising algebraic restrictions that occur.  \footnote{The author acknowledges support from U.S. National Science Foundation grants DMS 1107452, 1107263, 1107367 "RNMS: GEometric structures And Representation varieties" (the GEAR Network).}
\end{abstract}

\section*{Introduction}

\renewcommand\thesubsection{\arabic{subsection}}

Let $X$ be a Riemannian symmetric space of non-compact type, or a locally finite strongly transitive Euclidean building, with (type-preserving) isometry group $G$ and Weyl group $W$. Denote by $\ov{X}^g = X \cup \bord X$ the geodesic compactification of $X$ obtained by adding to $X$ its visual boundary $\bord X$, the space of classes of asymptotic geodesic rays (see~\cite[Chapter~II.8]{bridson_haefliger}). This article focuses on limits in $\ov{X}^g$ of convex subsets of $X$ (or more generally of $\CAT(0)$ spaces), as we shall see below.

\bigskip

Geometric limits have been studied by Harvey in the case of Fuchsian groups (see~\cite{harvey_discrete} and \cite{ceg}), using the Chabauty topology on the space of closed subgroups of a locally compact topological group (see~\cite{chabauty}). Denote by $\Flats(X)$ the space of all maximal flats in $X$, endowed with the topology induced by the Chabauty topology on the space ${\cal C}(X)$ of closed subsets of $X$. It is a homogeneous space under the action of $G$ on $X$: more specifically when $X$ is a symmetric space, the $G$-space $\Flats(X)$ is isomorphic to the homogeneous space $G/N_G(A)$, where $A$ is a maximal connected split semisimple abelian subgroup of $G$ and $N_G(A)$ is the normalizer of $A$ in $G$. For instance if $G=\PGL(n,\RRR)$, then $G/N_G(A)$ is the space of $n$-tuples in general position in $\RRR\PPP^{n-1}$. Many compactifications of such spaces have been defined, from an algebraic geometry point of view: the Fulton-McPherson compactification (see~\cite{fulton_macpherson}), the variety of reductions studied by Iliev, Manivel and Le Barbier Grünewald in the complex case (see~\cite{iliev_manivel_severi}, \cite{iliev_manivel}, \cite{lebarbier} and \cite{lebarbier_ex})... In another article, inspired by the work of Guivarc'h, Ji and Taylor (see~\cite{guivarch} and \cite{haettel_m2}), we have defined in the real case the Chabauty compactification of $G/N_G(A)$ inside the space of closed subgroups of $G$, and we have studied it when $G$ has real rank one, and for $G=\SL(3,\RRR)$ and $\SL(4,\RRR)$ (see~\cite{haettel_chabauty_cartan}).

\bigskip

The geometric and dynamic study of group actions on non-Riemannian homogeneous spaces $G/H$ (with non-compact isotropy group $H$) is growing (see for instance~\cite{OW}, \cite{kobayashi_exposition} and \cite{kassel_proper_actions}, even though the actions there are discontinuous). The case $H=A$ is very natural and has to be studied in detail; and given the crucial role of the asymptotic point of view when $H$ is compact, the study of the geometry at infinity of $G/A$ seems potentially fruitful. If $\Gamma$ is a lattice in a real semisimple Lie group $G$, studying the asymptotic geometry of $G/A$ could help to relate topological properties of $A$-orbits in $G/\Gamma$ to geometric properties of $\Gamma$-orbits in $G/A$: for instance, the study of limit sets of such $\Gamma$-orbits in the geometric compactification of $G/A$ defined in this article could bring information about the dual $A$-orbits in $G/\Gamma$.

\bigskip

Define the geometric compactification $\ov{\Flats(X)}^g$ of $\Flats(X)$ to be its closure inside the space ${\cal C}(\ov{X}^g)$ of closed subsets of $\ov{X}^g$, endowed with the Chabauty topology. If we have a divergent sequence of maximal flats of $X$, its limit in $\ov{\Flats(X)}^g$ represents what is asymptotically "seen" from a basepoint in $X$, in other words its visual limit (see Section~\ref{subsec:geometric_compactification_flats} for details).

Consider the structure of compact topological spherical building on $\bord X$, and denote by $\Cham(\bord X)$ the compact space of Weyl chambers of $\bord X$. Denote by $\Ap(X)$ the space of all apartments in $\bord X$. For each apartment in $\bord X$, its chambers are naturally indexed by the Weyl group $W$ of $X$, up to the $W$-action. Endow $\Ap(X)$ with the direct product topology $\Cham(\bord X)^W$, quotiented by the action of $W$. Define the geometric compactification $\ov{\Ap(X)}^g$ of $\Ap(X)$ to be its closure inside the compact space $\Cham(\bord X)^W$.

\bthmi \label{thmi:isomorphism} The natural $G$-equivariant homeomorphism between $\Flats(X)$ and $\Ap(X)$ extends to a $G$-equivariant homeomorphism between $\ov{\Flats(X)}^g$ and $\ov{\Ap(X)}^g$. \ethmi

Hence, in order to describe geometric limits of maximal flats in $X$, we only need to understand the geometric compactification of the space of apartments in the building at infinity of $X$, which is more combinatorial and more tractable.

Let ${\cal I}$ be a compact topological spherical building, let ${\cal A}$ be a fixed apartment of ${\cal I}$ and let $W$ be the Weyl group of ${\cal I}$. Every apartment of ${\cal I}$ is the image of ${\cal A}$ under a type-preserving, injective morphism of simplicial complexes from ${\cal A}$ to ${\cal I}$, hence we will in fact consider the space $\Mor_{inj}({\cal A},{\cal I})$ of such marked apartments of ${\cal I}$, which is a $W$-principal bundle over $\Ap({\cal I})$. Define the geometric compactification $\ov{\Mor_{inj}({\cal A},{\cal I})}^g$ of $\Mor_{inj}({\cal A},{\cal I})$ to be its closure inside the space $\Mor({\cal A},{\cal I})$ of (non-necessarily injective) morphisms from ${\cal A}$ to ${\cal I}$, endowed with the compact topology induced by the direct product space $\Cham({\cal I})^{\Cham({\cal A})}$.

If ${\cal I}$ is the join of two spherical buildings ${\cal I}_1$ and ${\cal I}_2$, then the geometric compactification of the space of marked apartments of ${\cal I}$ is just the Cartesian product of the geometric compactifications of the spaces of marked apartments of ${\cal I}_1$ and ${\cal I}_2$. Hence we only need to study the case where ${\cal I}$ is an irreducible building. We are able to describe $\ov{\Mor_{inj}({\cal A},{\cal I})}^g$ entirely in the following cases.

\bthmi \label{thmi:a1a2} Let ${\cal I}$ be of type $A_1$ or $A_2$, then $\ov{\Mor_{inj}({\cal A},{\cal I})}^g = \Mor({\cal A},{\cal I})$. \ethmi

This implies for instance the existence of visual limits of maximal flats in the symmetric space of unit volume ellipsoids in $\RRR^3$ which are not contained in any apartment at infinity.

\bigskip

If ${\cal I}$ is a spherical building of type $C_2$, call quadripod any union of four pairwise distinct Weyl chambers whose intersection is a vertex.

\bthmi \label{thmi:c2} Let $\KKK$ be a local field of characteristic different from $2$ (or a quaternion algebra over such a local field), let $V$ be a finite-dimensional (right) $\KKK$-vector space of dimension at least $5$. Consider a (possibly trivial) involutive automorphism $\sigma$ of $\KKK$, and let $q$ be a nondegenerate Hermitian form with respect to $\sigma$ on $V$ of Witt index $2$. Let ${\cal I}$ be the flag complex of totally isotropic subspaces of $V$: it is a classical thick spherical building of type $C_2$. Then every element of $\Mor({\cal A},{\cal I})$ whose image is not a quadripod belongs to $\ov{\Mor_{inj}({\cal A},{\cal I})}^g$. \ethmi

We give a simple cross-ratio condition characterizing which quadripods belong to the compactification $\ov{\Mor_{inj}({\cal A},{\cal I})}^g$. One should remark that the answer does not only depend on the Weyl group, but also on the local field $\KKK$. For instance, consider the case $\KKK=\RRR$ or $\CCC$, $V=\KKK^5$ with canonical basis $(e_1,\ldots,e_5)$ and $q(x)=x_1x_5+x_2x_4+x_3^2$. Consider four Weyl chambers $(\ell_i \subset \<e_1,e_2\>)_{1 \leq i \leq 4}$, where for all $i \in \lb 1,4\rb$, $\ell_i$ is a line in the plane $\<e_1,e_2\>$. The intersection of these Weyl chambers is the isotropic plane $\<e_1,e_2\>$, so their union is a quadripod called a plane-type quadripod. Then if $\KKK=\CCC$ it is always a limit of apartments, whereas if $\KKK=\RRR$ it is a limit of apartments if and only if the cross-ratio of $(\ell_1,\ell_2,\ell_3,\ell_4)$ is less than or equal to $1$.
%

\bigskip

Let $\KKK$ be a field, and let ${\cal I}$ be the spherical building of complete flags of the vector space $\KKK^4$, then an apartment of ${\cal I}$ is simply a tetrahedron in general position in $\PPP^3(\KKK)$. And a morphism from an apartment to ${\cal I}$ is given by $4$ points, $6$ lines and $4$ planes in $\PPP^3(\KKK)$ that satisfy the incidence conditions of a tetrahedron (see Section~\ref{sec:geometric_limits}).
We call a morphism of type $(L)$ if the $6$ lines are equal, and we call it of type $(XP)$ if the four points are equal and the four planes are equal.

\bthmi \label{thmi:a3} Let $\KKK$ be a local field, and let ${\cal I}$ be the spherical topological building of complete flags of the vector space $\KKK^4$. Then every element of $\Mor({\cal A},{\cal I})$ that is not of type $(L)$ nor $(XP)$ belongs to $\ov{\Mor_{inj}({\cal A},{\cal I})}^g$.

\smallskip

A element of $\Mor({\cal A},{\cal I})$ of type $(L)$ in general position belongs to $\ov{\Mor_{inj}({\cal A},{\cal I})}^g$ if and only if the cross-ratio of the four points is equal to the cross-ratio of the four planes.

A element of $\Mor({\cal A},{\cal I})$ of type $(XP)$ in general position belongs to $\ov{\Mor_{inj}({\cal A},{\cal I})}^g$ if and only if there is a projective involution of $\PPP^3(\KKK)$ fixing the common plane and the common point, and exchanging each line with its opposite.
\ethmi

In the Archimedean case, here is a corollary of these results in terms of visual limits of maximal flats. Let $G$ be a connected real semisimple Lie group with finite center, whose simple factors are of $\RRR$-rank one or locally isomorphic to $\SL(3,\RRR)$, $\SL(3,\CCC)$, $\SL(3,\HHH) \simeq \SU^*(6)$, $E_{6(-26)}$, $\SO(2,n)$, $\SU(2,n)$ or $\Sp(2,n)$ (where $n \geq 3$), $\SO(5,\CCC)$, $\SO^*(10)$, $\SL(4,\RRR)$ or $\SL(4,\CCC)$. Let $X$ be the symmetric space of $G$. Then we describe all the possible visual limits of divergent sequences of maximal flats in $X$, in terms of cross-ratios (see~Section~\ref{sec:geometric_limits} for precise statements).

\medskip

In the first section of the paper, we define the geometric compactifications of the spaces of flats and of apartments, and we prove Theorem~\ref{thmi:isomorphism}, by proving a general result on visual limits in $\CAT(0)$ spaces, Theorem~\ref{thm:limite_fermés}, which could be used to describe visual limits in other settings.

In the second section of the paper, we compute explicitly the geometric compactifications of the space of marked apartments in each of the cases of Theorems~\ref{thmi:a1a2}, \ref{thmi:c2} and \ref{thmi:a3}, and we emphasize the remarkably rich behaviours of visual limits of maximal flats that already occur in small ranks, and the algebraic restrictions on their existence that occur.

\medskip

\noindent {\small {\it Acknowledgements: } I would like to thank very warmly Frédéric Paulin, for his constant interest and precious help. I would also like to thank Yves Benoist, François Guéritaud and Fanny Kassel for useful discussions. I thank the referees for their useful reports.}

\section{Geometric compactifications of spaces of flats and apartments}

We begin by recalling the definition of a topological spherical building. Then we define in similar ways firstly the geometric compactification of the space of (marked) maximal flats of a symmetric space of non-compact type or a locally finite Euclidean building, and secondly the geometric compactification of the space of (marked) apartments of a compact topological spherical building. Then we show that the two compactifications are equivariantly isomorphic (Theorem~\ref{thm:compactification_plats}), by proving a more general result on visual limits in $\CAT(0)$ spaces (Theorem~\ref{thm:limite_fermés}).

Recall that if $E$ is a locally compact topological space, the set ${\cal C}(E)$ of closed subsets of $E$ is endowed with a natural compact topology called the Chabauty topology. For instance when $E$ is metrisable, then the Chabauty topology is metrizable for the pointed Gromov-Hausdorff distance: closed subsets are close if they almost coincide on large compact subsets. The Chabauty topology is also called Hausdorff or Vietoris topology, see~\cite{chabauty}, \cite{harpe_chabauty}, \cite[$\S 5$]{bourbaki_integration8}, \cite[Proposition~I.3.1.2, p.~59]{ceg}, \cite[Proposition~1.7, p.~58]{cdp}, \cite{RxZ}, \cite{haettel_m2} or \cite{haettel_chabauty_cartan}.

If $G$ is a locally compact topological group and $X$ is a locally compact topological  $G$-space, a $G$-compactification of $X$ is a pair $(\iota, K)$, where $K$ is a compact $G$-space and $\iota : X \hookrightarrow K$ is a $G$-equivariant topological embedding whith open and dense image.

By spherical or Euclidean building, we mean its spherical or Euclidean geometric realisation, and we consider its maximal apartments system (see~\cite{bridson_haefliger} and \cite{abramenko_brown}).

\subsection{Topological spherical buildings}

\label{subsec:top_buildings}

Given $k$ in $\NNN$ and $C$ a cellular complex, we denote by $C^{\langle k \rangle}$ the set of cells of $C$ with dimension $k$.

A \df{topological spherical building} (see~\cite{burns_spatzier} and \cite{ji_buildings}) is a spherical building ${\cal I}$ whose set of vertices ${\cal I}^{\langle 0 \rangle}$ is endowed with a topology such that for all $k \in \lb 0,d \rb$ (where $d$ is the dimension of ${\cal I}$), the set ${\cal I}^{\langle k \rangle}$ of $k$-simplices is closed for the topology induced by $\left({\cal I}^{\langle 0 \rangle}\right)^{k+1}$. If the space ${\cal I}^{\langle 0 \rangle}$ is (locally) compact, we say that the topological building ${\cal I}$ is (locally) compact.

Denote by ${\cal T}_{dist}$ the topology on ${\cal I}$ induced by the simplicial distance, which is not locally compact in general. The topology on the space $\Cham({\cal I}) = {\cal I}^{\langle d \rangle}$ of Weyl chambers of ${\cal I}$ defines a new topology ${\cal T}_{lc}$ on ${\cal I}$, coarser than ${\cal T}_{dist}$, for which a sequence $(x_n)_{n \in \NNN}$ converges to $x$ if and only if
\ben
\item there exists a sequence $(C_n)_{n \in \NNN}$ of Weyl chambers such that, for all $n \in \NNN$, we have $x_n \in C_n$, and such that the sequence $(C_n)_{n \in \NNN}$ converges to a Weyl chamber $C$ containing $x$ in $\Cham({\cal I})$;
\item if we denote by $\phi_n$ the type-preserving isometry from $C_n$ to $C$ for all $n \in \NNN$, then the sequence $(\phi_n(x_n))_{n \in \NNN}$ converges to $x$ in $C$.
\een
If the space ${\cal I}^{\langle 0 \rangle}$ is (locally) compact, then the topology ${\cal T}_{lc}$ is (locally) compact. An \df{automorphism of the topological building ${\cal I}$} is a (type-preserving) automorphism of the building ${\cal I}$ which is also a homeomorphism for the topology ${\cal T}_{lc}$.


\bigskip

For instance, let $\KKK$ be a \df{local field} (i.e. a field with a non-trivial valuation, complete and locally compact, thus isomorphic to $\RRR$, $\CCC$, a finite extension of $\QQQ_p$ or $\FFF_q((t))$, with $p$ prime and $q$ a power of a prime, see for instance~\cite[Chapitre~II]{serre_corpslocaux}), and let $G$ be the group of $\KKK$-points of an algebraic linear connected reductive $\KKK$-group $\underline{G}$. Let ${\cal I}$ be the spherical building of $(\underline{G},\KKK)$, then $G$ acts transitively on the vertices of ${\cal I}^{\langle 0 \rangle}$ of fixed type, with stabilizer a maximal proper parabolic $\KKK$-subgroup. Hence the set ${\cal I}^{\langle 0 \rangle}$ identifies with the disjoint union $\coprod_i G/P_i$, where the $P_i$'s are the groups of $\KKK$-points of maximal proper parabolic $\KKK$-subgroups containing a fixed minimal parabolic $\KKK$-subgroup. Endow ${\cal I}^{\langle 0 \rangle}$ with the topology induced by this disjoint union of $\KKK$-points of projective varieties: it is a compact space. This defines a structure of compact topological spherical building on ${\cal I}$, and the topology ${\cal T}_{lc}$ is metrisable.
%

%
\bigskip

In our cases, the topological spherical buildings will arise as visual boundaries of some $\CAT(0)$ spaces. Let $X$ be a product of symmetric spaces of non-compact type and locally finite strongly transitive Euclidean buildings. Then it is a complete locally compact $\CAT(0)$ space (see~\cite[Theorem~11.16, p.555]{abramenko_brown} for the building case). Its visual (or geodesic, or $\CAT(0)$) boundary at infinity $\bord X$ has a natural structure ${\cal I}$ of spherical building (see~\cite{bgs} and \cite[§~11.8]{abramenko_brown}). The space ${\cal I}^{\langle 0 \rangle}$ of vertices of ${\cal I}$ is included in $\bord X$, let us endow it with the induced compact topology. This defines a natural structure of compact topological spherical building on ${\cal I}$, and the topology ${\cal T}_{lc}$ on ${\cal I}$ is the topology of $\bord X$. Furthermore ${\cal I}$ is topologically strongly transitive, which means that its group of (topological) automorphisms acts transitively on the set of flags $(C \subset {\cal A})$, where $C$ is a Weyl chamber and ${\cal A}$ is an apartment of ${\cal I}$.

\subsection{Geometric compactification of the space of (marked) maximal flats}

\label{subsec:geometric_compactification_flats}

Let $X$ be a product of symmetric spaces of non-compact type and of locally finite strongly transitive Euclidean buildings. Let $G$ be the product of the isometry groups of the symmetric space factors and of the type-preserving automorphism groups of the building factors. Denote by $\Flats(X)$ the set of all maximal flats of $X$ (i.e. the set of apartments if $X$ is a building), endowed with the topology induced by the Chabauty topology on the space ${\cal C}(X)$ of closed subsets of $X$, and with the $G$-action on the left. Denote by $W$ the Weyl group of $X$.

Consider the geodesic $G$-compactification $\ov{X}^g = X \cup \bord X$ of $X$. Consider the $G$-equivariant embedding
\beq \Flats(X) & \lra & {\cal C}(\ov{X}^g) \\
F & \lmapsto & \ov{F},\eeq
where $\ov{F}$ denotes the topological closure of $F$ in $\ov{X}^g$. Call \df{geometric compactification} of the space $\Flats(X)$ of maximal flats of $X$ the closure $\ov{\Flats(X)}^g$ of its image in ${\cal C}(\ov{X}^g)$.

\bigskip

We can define as well a geometric compactification of the space of marked flats of $X$. Denote by $\Cham(\bord X)$ the space of all closed Weyl chambers in the spherical building at infinity $\bord X$, endowed with the topology induced by the Chabauty topology on the space ${\cal C}(\bord X)$ of closed subsets of $\bord X$. Call \df{marked flat} of $X$ any $(F,C) \in \Flats(X) \times \Cham(\bord X)$ such that $C \subset \bord F$. Denote by $\Flats_m(X)$ the set of all marked flats of $X$, endowed with the compact topology induced by the product topology on $\Flats(X) \times \Cham(\bord X)$. It has a natural $(G \times W)$-action on the left. Consider the $G$-equivariant embedding
\beq \Flats_m(X) & \lra & {\cal C}(\ov{X}^g) \times \Cham(\bord X)\\
(F,C) & \lmapsto & (\ov{F},C),\eeq
and call \df{geometric compactification} of the space $\Flats_m(X)$ of marked flats of $X$ the closure $\ov{\Flats_m(X)}^g$ of its image.

\bigskip

The forgetful map
\beq \Flats_m(X) & \lra & \Flats(X) \\
(F,C) & \lmapsto & F\eeq
is continuous, surjective, $G$-equivariant and with finite fibers. It is the quotient map under the action of $W$.

\bigskip

If $X = \prod_{i=1}^k X_i$ is a product, where each $X_i$ is a product of symmetric spaces of non-compact type and locally finite strongly transitive Euclidean buildings, then there is a natural $(G \times W)$-equivariant homeomorphism
\beq \prod_{i=1}^k \Flats_m(X_i) & \lra & \Flats_m(X) \\
(F_i,C_i)_{1 \leq i \leq k} & \lmapsto & (\prod_{i=1}^k F_i,\star_{i=1}^k C_i),\eeq
where $\star_{i=1}^k C_i$ denotes the simplicial join of $(C_1,\ldots,C_k)$, whose set of (possibly empty) simplices is the product of the set of (possibly empty) simplices of each $C_i$ (see for instance~\cite[Definition~I.7A.2]{bridson_haefliger}).

This extends to a $(G \times W)$-equivariant homeomorphism between $\prod_{i=1}^k \ov{\Flats_m(X_i)}^g$ and $\ov{\Flats_m(X)}^g$, and similarly between $\prod_{i=1}^k \ov{\Flats(X_i)}^g$ and $\ov{\Flats(X)}^g$. Hence we only need to study the case where $X$ is irreducible.

\subsection{Geometric compactification of the space of (marked) apartments}

Let ${\cal I}$ be a compact topological spherical building (see Subsection~\ref{subsec:top_buildings}) of dimension $d$ and topological automorphism group $G$. Fix an apartment ${\cal A}$ of ${\cal I}$, and denote by $W$ its Weyl group. Recall that a morphism from ${\cal A}$ to ${\cal I}$ is a (non-necessarily injective) type-preserving morphism of typed simplicial complexes. Denote by $\Mor({\cal A},{\cal I})$ the set of all morphisms from ${\cal A}$ to ${\cal I}$. It has a natural $(G \times W)$-action on the left, given by
$$ \forall g \in G, \forall w \in W, \forall f \in \Mor({\cal A},{\cal I}), (g,w) \cdot f = g \circ f \circ w^{-1}.$$

Call \df{marked apartment} of ${\cal I}$ any morphism from ${\cal A}$ to ${\cal I}$ whose image is an apartment of ${\cal I}$. Denote by $\Mor_{inj}({\cal A},{\cal I})$ or $\Ap_m({\cal I})$ the subset of all marked apartments of ${\cal I}$. This notation is justified by the fact that a morphism from ${\cal A}$ to ${\cal I}$ is injective if and only if its image is an apartment (see~\cite[Proposition~4.59, p.193]{abramenko_brown}).

A morphism from ${\cal A}$ to ${\cal I}$ is characterized by the images of the finite number of Weyl chambers $\Cham({\cal A})$ of ${\cal A}$. Denote by $\Cham({\cal I}) = {\cal I}^{\langle d \rangle}$ the space of closed Weyl chambers of ${\cal I}$, endowed with the compact topology induced by the product topology on $({\cal I}^{\langle 0 \rangle})^{d+1}$.

Endow $\Mor({\cal A},{\cal I})$ with the topology  induced by the product topology on $\Cham({\cal I})^{\Cham({\cal A})}$, which is the same as the topology induced by the compact-open topology for the topology ${\cal T}_{lc}$ on ${\cal I}$. An element $f$ of $\Cham({\cal I})^{\Cham({\cal A})}$ is a morphism if and only if for all Weyl chambers $C,C'$ of ${\cal A}$ whose intersection is a non-empty facet, $f(C) \cap f(C')$ is a non-empty facet of the same type. Those are a finite number of closed conditions, so $\Mor({\cal A},{\cal I})$ is compact.

Call \df{geometric compactification} of the space $\Mor_{inj}({\cal A},{\cal I})=\Ap_m({\cal I})$ of marked apartments of ${\cal I}$ its closure $\ov{\Mor_{inj}({\cal A},{\cal I})}^g=\ov{\Ap_m({\cal I})}^g$ inside $\Mor({\cal A},{\cal I})$. Denote by $\partial \Ap_m({\cal I}) = \ov{\Ap_m({\cal I})}^g \bs \Ap_m({\cal I})$ its boundary.

\bigskip

We can define as well a geometric compactification of the space $\Ap({\cal I})$ of all (unmarked) apartments of ${\cal I}$, endowed with the quotient topology of $\Ap_m({\cal I})$, which is the same as the compact topology induced by the Chabauty topology on the space ${\cal C}({\cal I},{\cal T}_{lc})$ of closed subsets of ${\cal I}$ endowed with the topology ${\cal T}_{lc}$. It has a natural $G$-action on the left. Call \df{geometric compactification} of the space $\Ap({\cal I})$ of apartments of ${\cal I}$ its closure $\ov{\Ap({\cal I})}^g$ inside ${\cal C}({\cal I},{\cal T}_{lc})$. Denote by $\partial \Ap({\cal I}) = \ov{\Ap({\cal I})}^g \bs \Ap({\cal I})$ its boundary.

\bigskip

The forgetful map
\beq \Ap_m({\cal I}) & \lra & \Ap({\cal I}) \\
(f : {\cal A} \ra {\cal I}) & \lmapsto & f({\cal A})\eeq
is continuous, surjective, $G$-equivariant and with finite fibers. It is the quotient map under the action of $W$.

\bigskip

If ${\cal I} = \star_{i=1}^k {\cal I}_i$ is a join of compact topological spherical buildings, then there is a natural $(G \times W)$-equivariant homeomorphism
\beq \prod_{i=1}^k \Ap_m({\cal I}_i) & \lra & \Ap_m({\cal I}) \\
(f_i)_{1 \leq i \leq k} & \lmapsto & \star_{i=1}^k f_i,\eeq
where $\star_{i=1}^k f_i$ denotes the simplicial join of the simplicial maps $(f_i)_{1 \leq i \leq k}$.

This extends to a $(G \times W)$-equivariant homeomorphism between $\prod_{i=1}^k \ov{\Ap_m({\cal I}_i)}^g$ and $\ov{\Ap_m({\cal I})}^g$, and similarly between $\prod_{i=1}^k \ov{\Ap({\cal I}_i)}^g$ and $\ov{\Ap({\cal I})}^g$. Hence we only need to study the case where ${\cal I}$ is irreducible.

\bigskip

In the classical case, there is an algebraic interpretation of this compactification. Let $\KKK$ be a local field, and let $G$ be the group of $\KKK$-points of an algebraic linear connected reductive $\KKK$-group $\underline{G}$, with (restricted) Weyl group $W$. Let ${\cal I}$ be the compact topological spherical building of $(\underline{G},\KKK)$, and let $A$ be a maximal $\KKK$-split torus of $G$. Then the space $\Ap_m({\cal I})$ is $(G \times W)$-equivariantly homeomorphic to the homogeneous space $G/Z_G(A)$. Fix $P$ a minimal parabolic $\KKK$-subgroup of $G$ containing $A$. Then the embedding
\beq \Ap_m({\cal I}) \simeq G/Z_G(A) & \lra & \prod_{w \in W} G/wPw^{-1} \\
gZ_G(A) & \lmapsto & (gwPw^{-1})_{w \in W} \eeq
extends to a $(G \times W)$-equivariant embedding of $\ov{\Ap_m({\cal I})}^g$ into the algebraic projective $\KKK$-variety $\prod_{w \in W} G/wPw^{-1}$.

\subsection{Geometric limits of closed subsets in $\CAT(0)$ spaces}

Let $X$ be a complete, locally compact $\CAT(0)$ metric space, let $\bord X$ be its geodesic (or visual, or $\CAT(0)$) boundary and $\ov{X}^g=X \cup \bord X$ its geodesic compactification (see~\cite{bridson_haefliger} and \cite{bgs}). Assume there exists a non-empty set ${\cal F}$ of closed non-empty subsets of $X$ satisfying the following two properties.
\ben
\item \label{condition_cocompact} The subspace $\{(F,x) \,:\, F \in {\cal F},x \in F\}$ of pointed elements of ${\cal F}$ is closed in the space ${\cal C}(X) \times X$, endowed with the product topology, and it is $\Isom(X)$-cocompact (recall that ${\cal C}(X)$ denotes the space of closed subsets of $X$, endowed with the Chabauty topology).
\item \label{condition_oppose} For all $F \in {\cal F}$ and all $\xi \in \bord X$, there exists $\eta \in \bord F$ \df{opposite} to $\xi$, that is there exists a geodesic in $X$ whose endpoints are $\eta$ and $\xi$.
\een

\bigskip

Here are families of examples where these two properties are satisfied. In each case ${\cal F}$ is the set of all maximal flats in $X$ and the boundary at infinity $\bord X$ is naturally a spherical building, which implies the opposition condition~(\ref{condition_oppose}).
\bit
\item $X$ is a symmetric space of non-compact type.
\item $X$ is a locally finite strongly transitive Euclidean building.
\item $X$ is a Gromov hyperbolic complete locally compact $\CAT(0)$ metric space, with extendible geodesics, whose isometry group acts cocompactly on pointed geodesics in $X$.
\eit

When $\bord X$ is a spherical building, we can also consider for ${\cal F}$ a set of closed subsets of $X$ containing a maximal flat, which also implies the opposition condition~(\ref{condition_oppose}).
\bit
\item $X$ is a Hermitian symmetric space of non-compact type, with ${\cal F}$ the set of all maximal polydiscs in $X$.
\item $X$ is a locally finite hyperbolic building (see~\cite{bourdon}, \cite{gaboriau_paulin}), whose isometry group acts strongly transitively, with ${\cal F}$ the set of all apartments of $X$.
\eit

Furthermore, any product of finitely many of the examples above (with the $\ell^2$ product metric) satisfies the two properties.

\bigskip

Let us denote by $\sa_x$ the visual angle on $\ov{X}^g$ at $x \in X$ (see for instance~\cite[Chapter~II.3]{bridson_haefliger}).

\blem \label{lem:angle_zero} Fix $x \in X$. Let $(y_n)_{n \in \NNN}$ be a sequence in $X$ converging to $\xi \in \bord X$, and let $\eta \in \bord X$ be opposed to $\xi$. Then the sequence $(\sa_{y_n}(x,\eta))_{n \in \NNN}$ converges to $0$. \elem

\begin{figure}[!h]
\begin{center}
\includegraphics[height=6cm]{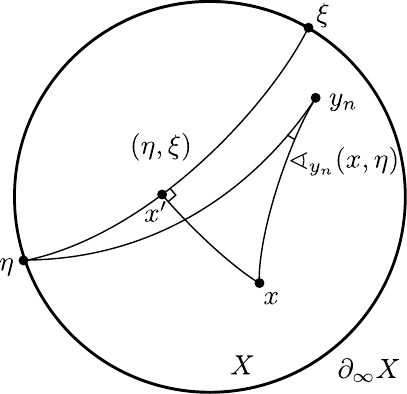}
\caption{Lemma~\ref{lem:angle_zero}.}
\label{fig:lemme_angles}
\end{center}
\end{figure}

\bp Fix a geodesic $(\eta, \xi)$ in $X$ whose endpoints are $\eta$ and $\xi$, and denote by $x'$ the orthogonal projection of $x$ to $(\eta, \xi)$ (see Figure~\ref{fig:lemme_angles}). We know that the sequence $(\sa_{x'}(y_n,\eta))_{n \in \NNN}$ converges to $\sa_{x'}(\xi,\eta)=\pi$.

Since for all $n \in \NNN$ we have $\sa_{x'}(y_n,\eta)+\sa_{y_n}(x',\eta) \leq \pi$, we deduce that the sequence $(\sa_{y_n}(x',\eta))_{n \in \NNN}$ converges to $0$.

As $(\sa_{y_n}(x,x'))_{n \in \NNN}$ converges to $0$, we conclude that the sequence $(\sa_{y_n}(x,\eta))_{n \in \NNN}$ converges to $0$. \ep

\bthm \label{thm:limite_fermés} Let $(F_n)_{n \in \NNN}$ be a sequence in ${\cal F}$ leaving every compact subset of $X$, and such that the sequence $(\bord F_n)_{n \in \NNN}$ converges to a closed subset $C$ of $\bord X$ in $\mathcal{C}(\bord X)$. Then the sequence $(\ov{F_n})_{n \in \NNN}$ converges to $C$ in $\mathcal{C}(\ov{X}^g)$. \ethm

\begin{figure}[!h]
\begin{center}
\includegraphics[height=6cm]{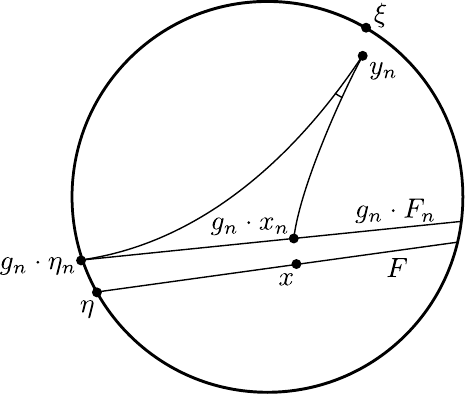}
\caption{Theorem~\ref{thm:limite_fermés}.}
\label{fig:theoreme_angles}
\end{center}
\end{figure}

\bp Since the sequence $(F_n)_{n \in \NNN}$ leaves every compact subset of $X$, every accumulation point of the sequence $(\ov{F_n})_{n \in \NNN}$ is included in $\bord X$. Let $(x_n)_{n \in \NNN}$ be a sequence of points of $(F_n)_{n \in \NNN}$ converging to $\zeta \in \bord X$. Since $\Isom(X)$ acts cocompactly on pointed elements of ${\cal F}$, we may assume by the cocompactness condition~(\ref{condition_cocompact}) that there exists $F \in {\cal F}$, $x \in F$ and a sequence $(g_n)_{n \in \NNN}$ in $\Isom(X)$ such that the sequence $(g_n \cdot F_n)_{n \in \NNN}$ converges to $F$ in ${\cal C}(X)$ and that the sequence $(g_n \cdot x_n)_{n \in \NNN}$ converges to $x$ in $X$ (see Figure~\ref{fig:theoreme_angles}).

Up to passing to a subsequence, we may assume that the sequence $(y_n = g_n \cdot x)_{n \in \NNN}$ converges to a point $\xi \in \bord X$, as this sequence leaves every compact of $X$:
$$d(x,y_n) = d(g_n^{-1} \cdot x,x) \geq d(x_n,x) - d(g_n^{-1} \cdot x,x_n) \geq d(F_n,x) - d(x,g_n \cdot x_n) \ral{n \ra +\infty} +\infty.$$

The opposition condition~(\ref{condition_oppose}) ensures the existence of $\eta \in \bord F$ opposed to $\xi$. We apply Lemma~\ref{lem:angle_zero} to the basepoint $x \in X$, the sequence $(y_n)_{n \in \NNN}$ and the point $\eta \in \bord F$: the sequence $(\sa_{y_n}(x,\eta))_{n \in \NNN}$ converges to $0$.

Since the sequence $(g_n \cdot F_n)_{n \in \NNN}$ converges to $F$, choose for all $n \in \NNN$ an element $\eta_n \in \bord F_n$ such that the sequence $(g_n \cdot \eta_n)_{n \in \NNN}$ converges to $\eta$ in $\bord X$. Then we have
$$\sa_{y_n}(g_n \cdot x_n,g_n \cdot \eta_n) \leq \sa_{y_n}(g_n \cdot x_n,x)+\sa_{y_n}(x,\eta)+\sa_{y_n}(\eta,g_n \cdot \eta_n) \ral{n \ra +\infty} 0.$$

Since $g_n^{-1}$ is an isometry of $X$ for all $n \in \NNN$, we deduce that the sequence $(\sa_{x}(x_n,\eta_n))_{n \in \NNN}$ converges to $0$. As the sequence $(x_n)_{n \in \NNN}$ converges to $\zeta$, we conclude that the sequence $(\sa_{x}(\zeta,\eta_n))_{n \in \NNN}$ converges to $0$, hence the sequence $(\eta_n)_{n \in \NNN}$ converges to $\zeta$ in $\bord X$. So $\zeta \in C$. The other inclusion being obvious, we conclude that the sequence $(\ov{F_n})_{n \in \NNN}$ converges to $C$ in $\mathcal{C}(\ov{X}^g)$.\ep

\subsection{The isomorphism between the compactifications}

Let $(X,G)$ be as in the beginning of Subsection~\ref{subsec:geometric_compactification_flats}. Denote by ${\cal I} = \bord X$ the compact topological spherical building at infinity of $X$. Denote by $d$ the dimension of ${\cal I}$. It is equal to $r-1$, where $r$ is the real rank of $X$ if $X$ is a symmetric space, the dimension of $X$ minus $1$ if $X$ is a Euclidean building, or their sum if $X$ is a product. Fix a maximal flat $F_0$ in $X$, denote by ${\cal A}=\bord F_0$ its apartment at infinity, and fix a Weyl chamber $C_0$ of ${\cal A}$.

\medskip

Consider the natural $G$-equivariant and $(G \times W)$-equivariant homeomorphisms
\beq \iota : \Flats(X) & \lra & \Ap({\cal I}) \\
F & \lmapsto & \bord F\\
\mbox{and } \iota_m : \Flats_m(X) & \lra & \Ap_m({\cal I}) \\
(F,C) & \lmapsto & f\in \Mor_{inj}({\cal A},{\cal I}) \mbox{ such that } f({\cal A})=\bord F \mbox{ and } f(C_0)=C,\eeq

Consider furthermore the map $\phi : \ov{\Ap({\cal I})}^g \ra {\cal C}(\ov{X}^g)$ defined as follows. If $A \in \Ap({\cal I})$, then define $\phi(A)$ to be $\ov{F}$, where $F$ is the unique maximal flat in $X$ such that $\bord F= A$. And if $A \in \partial \Ap({\cal I})$, define $\phi(A) = A \subset {\cal I}$.

Consider as well the marked version of $\phi$, that is the map $\phi_m : \ov{\Ap_m({\cal I})}^g \ra {\cal C}(\ov{X}^g) \times \Cham({\cal I})$ defined as follows. If $f \in \Ap_m({\cal I})$, then define $\phi_m(f)$ to be $(\ov{F},f(C_0))$ where $F$ is the unique maximal flat in $X$ such that $\bord F = f({\cal A})$. And if $f \in \partial \Ap_m({\cal I})$, define $\phi_m(f)$ to be $(f({\cal A}),f(C_0))$.

\bthm \label{thm:compactification_plats} The geometric compactification $\ov{\Flats_m(X)}^g$ of the space of marked flats of $X$ is $(G \times W)$-isomorphic to the geometric compactification $\ov{\Ap_m({\cal I})}^g$ of the space of marked apartments of ${\cal I}$. More precisely, the following diagram of $(G \times W)$-equivariant embeddings is well-defined and commutes.
\begin{equation} \begin{array}{ccc} \Flats_m(X) & \lhookra & \ov{\Flats_m(X)}^g \\
\iota_m \downarrow \! \wr &  & \wr \! \uparrow \phi_m \\
\Ap_m({\cal I}) & \lhookra & \ov{\Ap_m({\cal I})}^g.\end{array}
\label{eqn:diagramme}
\end{equation}
\ethm

Taking the quotient with respect to $W$ gives the same unmarked result.

\bcor \label{cor:compactification_plats} The geometric compactification $\ov{\Flats(X)}^g$ of the space of flats of $X$ is $G$-isomorphic to the geometric compactification $\ov{\Ap({\cal I})}^g$ of the space of apartments of ${\cal I}$. More precisely, the following diagram of $G$-equivariant embeddings is well-defined and commutes.
$$ \begin{array}{cccr} \Flats(X) & \lhookra & \ov{\Flats(X)}^g & \\
\iota \downarrow \! \wr &  & \wr \! \uparrow \phi & \\
\Ap({\cal I}) & \lhookra & \ov{\Ap({\cal I})}^g. &  \end{array}$$ \end{cor}

\bp It is clear that the diagram~(\ref{eqn:diagramme}) commutes, and that each arrow is $(G \times W)$-equivariant and injective. So it is enough to show that $\phi_m$ is continuous, with values in $\ov{\Flats_m(X)}^g$.

\bigskip

The map $\phi_m$ restricted to $\Ap_m({\cal I})$ is a homeomorphism onto $\Flats_m(X)$, since both spaces are endowed with the topology of $G$-homogeneous spaces. The topology on each space of closed Weyl facets ${\cal I}^{\langle k \rangle}$ is induced by the Chabauty topology on ${\cal C}(\bord X)$, hence $\phi_m$ restricted to $\partial \Ap_m({\cal I})$ is an embedding into ${\cal C}(\bord X)$.

\bigskip

Let us show that $\phi_m$ is continuous on $\ov{\Ap_m({\cal I})}^g$. Since $\ov{\Ap_m({\cal I})}^g$ and ${\cal C}(\ov{X}^g) \times \Cham({\cal I})$ are metrisable, we need only prove the sequential continuity. Let $f \in \partial \Mor_{inj}({\cal A},{\cal I})$, and let $(f_n)_{n \in \NNN}$ be a sequence in $\Mor_{inj}({\cal A},{\cal I})$ converging to $f$ in $\Mor({\cal A},{\cal I})$. For all $n \in \NNN$, let $(F_n,C_n)=\phi_m(f_n)$.

The sequence $(F_n)_{n \in \NNN}$ of maximal flats of $X$ goes to infinity in $\Flats(X)$, and the sequence $(\bord F_n)_{n \in \NNN}$ converges to $f({\cal A})$ in ${\cal C}(\bord X)$ by assumption. Applying Theorem~\ref{thm:limite_fermés}, we deduce that the sequence $(\ov{F_n})_{n \in \NNN}$ converges to $f({\cal A})$ in ${\cal C}(\ov{X}^g)$.

Furthermore by assumption we know that the sequence $(C_n=f_n(C_0))_{n \in \NNN}$ of Weyl chambers converges to $f(C_0)$. Hence the sequence $((F_n,C_n)=\phi_m(f_n))_{n \in \NNN}$ converges to $(f({\cal A}),f(C_0))=\phi_m(f)$ in ${\cal C}(\ov{X}^g) \times \Cham(\bord X)$. So $\phi_m$ is continuous at $f$, and furthermore $\phi_m(f) \in \ov{\Flats_m(X)}^g$.

\bigskip

The map $\phi_m$ is continuous, injective, from the compact space $\ov{\Ap_m({\cal I})}^g$ into the Hausdorff space $\ov{\Flats_m(X)}^g$ : it is an embedding. Its image is a compact subspace of $\ov{\Flats_m(X)}^g$, which contains the dense subspace $\phi_m(\Ap_m({\cal I})) = \Flats_m(X)$, hence $\phi_m$ is surjective. So $\phi_m$ is a $(G \times W)$-equivariant homeomorphism from $\ov{\Ap_m({\cal I})}^g$ onto $\ov{\Flats_m(X)}^g$.
\ep

\section{Geometric limits of marked apartments in topological spherical buildings}

In this section, we compute explicitly the geometric compactifications of the space of marked apartments in each of the cases of Theorems~\ref{thmi:a1a2}, \ref{thmi:c2} and \ref{thmi:a3}. First we look at the rank $1$ case. Then in the rank $2$ case (or the topologically Moufang polygon case), we study completely the type $A_2$, and for the type $C_2$ we get a complete description in the orthogonal and Hermitian classical cases. Finally, we study the case of $\PGL(4)$ over a field.

\label{sec:geometric_limits}

\subsection{Rank $1$}

A compact topological spherical building ${\cal I}$ of rank $1$ is just a compact space ${\cal I}^{\langle 0 \rangle}$ with the set of unordered pairs of distinct points of ${\cal I}^{\langle 0 \rangle}$ as apartment system. So the space of marked apartments of ${\cal I}$ is the set of ordered pairs of distinct points of ${\cal I}^{\langle 0 \rangle}$.

Assume that ${\cal I}^{\langle 0 \rangle}$ has no isolated point, then the geometric compactification $\ov{\Mor_{inj}({\cal A},{\cal I})}^g$ of the space of marked apartments of ${\cal I}$ is just $\Mor({\cal A},{\cal I}) = ({\cal I}^{\langle 0 \rangle})^2$. And the geometric compactification $\ov{\Ap({\cal I})}^g$ of the space of unmarked apartments of ${\cal I}$ is the quotient of $({\cal I}^{\langle 0 \rangle})^2$ by the diagonal involution. This is the first part of Theorem~\ref{thmi:a1a2}.

For instance, the space of oriented geodesics of the real hyperbolic plane $\HHH^2_\RRR$ is homeomorphic to the one-sheeted hyperboloid $\{(x,y,z) \in \RRR^3 \,:\, x^2+y^2-z^2=1\}$ (in the hyperboloid model of $\HHH^2_\RRR$), or to the space of ordered pairs of distinct points of $\SS^1$ (in the disc model), and the geometric compactification is homeomorphic to the torus $(\SS^1)^2$. And the geometric compactification of the space of non-oriented geodesics of $\HHH^2_\RRR$ is homeomorphic to the quotient of $(\SS^1)^2$ by the diagonal involution.

If we apply Theorem~\ref{thm:compactification_plats}, we get the following well-known result.

\bcor Let $X$ be a $\RRR$-rank $1$ symmetric space of non-compact type or a locally finite, strongly transitive tree. Then any visual limit of a divergent sequence of geodesics in $X$ is a single point in the visual boundary $\bord X$. \ecor

\subsection{Rank $2$ : topologically Moufang polygons}

\label{subsec:rank_2}

\subsubsection{Notation}

Let ${\cal I}$ be a topological spherical building of rank $2$ : it is a bipartite graph, with a topology on its set of edges. Assume that this topology is locally compact and has no isolated point. Let $G$ be the group of automorphisms of the topological spherical building ${\cal I}$, assume that ${\cal I}$ is topologically strongly transitive. Let ${\cal I}^{\langle 0 \rangle}$ be the space of vertices of ${\cal I}$, and let ${\cal I}^{\langle 1 \rangle}$ be the space of edges of ${\cal I}$.

Fix an apartment ${\cal A}$ of ${\cal I}$: it is the boundary of a $2p$-gon, with $p \in \NNN$ at least $2$. We will be interested in the cases $p=2$, $p=3$ and $p=4$. Let $(C_i)_{i \in \ZZZ/2p\ZZZ}$ be the edges (Weyl chambers) of ${\cal A}$, ordered in such a way that for all $i \in \ZZZ/2p\ZZZ$, the edges $C_i$ and $C_{i+1}$ are adjacent. For all $i \in \ZZZ/2p\ZZZ$, let $x_{i,i+1}$ be the intersection of $C_i$ and $C_{i+1}$ ; the parity of $i$ determines the type of $x_{i,i+1}$.


If $\alpha$ is a root (or half-apartment) of ${\cal I}$, let $U_\alpha$ be its root group :
$$ U_\alpha =\{ g \in G \,:\, g(\alpha) = \alpha\}.$$
The topological building ${\cal I}$ is said to be \df{topologically Moufang} if, for every root $\alpha$ of ${\cal I}$, the root group $U_\alpha$ acts transitively on the set of apartments containing $\alpha$. Note that since $U_\alpha$ is the topological root group, this definition is a priori stronger than the Moufang condition for a general spherical building. Note that Grundhöfer, Knarr, Kramer, Van Maldeghem and Weiss classified compact irreducible spherical buildings of rank at least $2$ which are Moufang (see~\cite{kramer_flag_homogeneous}, \cite{kramer_flag_homogeneous_II} and \cite{kramer_compact_moufang_buildings}), and in particular it follows that they are also topologically Moufang. This topologically Moufang condition is stable under product, and it is verified for the topological spherical building of a classical group over a local field (see~\cite[Chapter~7]{abramenko_brown}). In Subsection~\ref{subsec:rank_2}, we will assume that ${\cal I}$ is topologically Moufang.

\medskip

Since by assumption ${\cal I}^{\langle 1 \rangle}$ has no isolated point, each root group of ${\cal I}$ is non-compact.

\medskip

For all $i \in \ZZZ/2p\ZZZ$, let $\alpha_i$ be the root (or half-apartment) $\{C_i,C_{i+1},\ldots,C_{i+p-1}\}$, and let $U_i$ be its root group $U_{\alpha_i}$. The Weyl group $W$ of ${\cal I}$ is isomorphic to the dihedral group $D_{2p}$ of order $2p$.

If $x$ is a vertex in ${\cal I}$, let $\st(x)$ be the \df{star} of $x$:
$$ \st{x} = \{ C \in {\cal I}^{\langle 1 \rangle} \,:\, x \in C\},$$
it is a closed subset of ${\cal I}$ for ${\cal T}_{lc}$.

\bigskip

Fix a root $\alpha_i$ of ${\cal A}$, and consider the folding $p_i : {\cal A} \ra \alpha_i$ onto $\alpha_i$. The following lemma is immediate, but quite useful.

\blem \label{lem:proj} Let $(g_n)_{n \in \NNN}$ be a divergent sequence in $U_i$. Then the sequence $(g_n|_{\cal A})_{n \in \NNN}$ converges to $p_i$ in $\Mor({\cal A},{\cal I})$. Consequently, if $f \in \ov{\Mor_{inj}({\cal A},{\cal I})}^g$ and if $f({\cal A}) \subset {\cal A}$, then $p_i \circ f = \liml_{n \ra +\infty} g_n \circ f \in \ov{\Mor_{inj}({\cal A},{\cal I})}^g$ for all $i \in \ZZZ/2p\ZZZ$. \elem

%

\subsubsection{Type $A_1^2$ : product of topologically Moufang lines}

Assume here that $p=2$. Then the building ${\cal I}$ is not irreducible, so we can deduce the following proposition from the $A_1$ case, but its direct proof is easy.

\bpro \label{pro:limites_a12} If ${\cal I}$ has type $A_1^2$, then the space $\Mor({\cal A},{\cal I})$ has 4 orbits under the action of $(G \times W)$, whose representatives are
$$ \id : {\cal A} \ra {\cal A}, p_1 : {\cal A} \ra \alpha_1 , p_2 : {\cal A} \ra \alpha_2 \mbox{ and } p_1 \circ p_2 = p_2 \circ p_1 : {\cal A} \ra \alpha_1 \cap \alpha_2 = C_2.$$
According to Lemma~\ref{lem:proj}, we deduce that $\ov{\Mor_{inj}({\cal A},{\cal I})}^g = \Mor({\cal A},{\cal I})$.
\epro

If we apply Theorem~\ref{thm:compactification_plats}, we get the following result.

\bcor Let $X$ be a product of two $\RRR$-rank $1$ symmetric spaces of non-compact type or locally finite, strongly transitive trees. Then the visual limits of divergent sequences of maximal flats in $X$ are the two types of half-apartments at infinity, or a Weyl chamber in the visual boundary $\bord X$. \ecor

\subsubsection{Type $A_2$ : topologically Moufang planes}

Assume here that $p=3$. For instance, we may take $\KKK$ to be a finite-dimensional division algebra over a local field, and ${\cal I}$ the topological spherical building of complete flags of the right vector space $\KKK^3$.


\medskip

Call (marked) \df{tripod} of ${\cal I}$ any morphism $f : {\cal A} \ra {\cal I}$ whose image consists of three pairwise intersecting edges. Call \df{type} of a tripod the type of this common intersection.

For example, fix $C$ an edge in ${\cal I}$ containing $x_{12}$, different from $C_1$ and $C_2$. Then $t_{12} : {\cal A} \ra {\cal I}$, defined by $t_{12}(C_6)=t_{12}(C_1)=C_1$, $t_{12}(C_2)=t_{12}(C_3)=C_2$ and $t_{12}(C_4)=t_{12}(C_5)=C$, is a tripod of one type.

\begin{figure}[!h]
\begin{center}
\includegraphics[height=4cm]{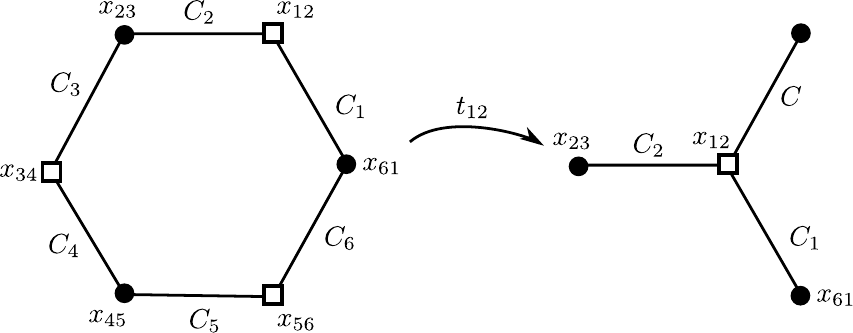}
\includegraphics[height=4cm]{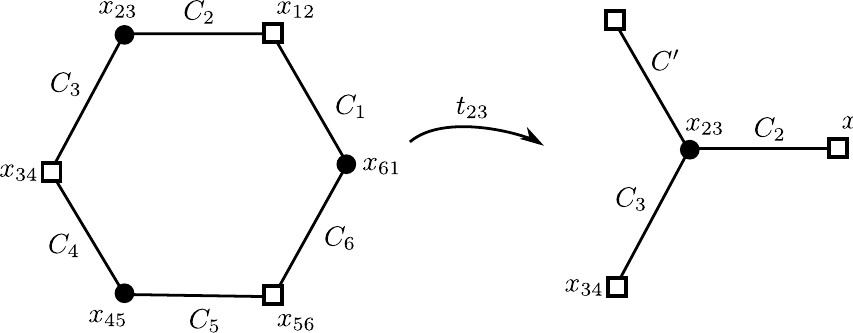}
\caption{The marked tripods of the two types}
\label{fig:marked_tripods}
\end{center}
\end{figure}

Similarly, fix $C'$ an edge in ${\cal I}$ containing $x_{23}$, different from $C_2$ and $C_3$. Then $t_{23} : {\cal A} \ra {\cal I}$, defined by $t_{23}(C_1)=t_{23}(C_2)=C_2$, $t_{23}(C_3)=t_{23}(C_4)=C_3$ and $t_{23}(C_5)=t_{23}(C_6)=C'$, is a tripod of the other type (see Figure~\ref{fig:marked_tripods}).

Two elements $f$, $f'$ in $\Mor({\cal A},{\cal I})$ are said to be \emph{equivalent} if there exist $w \in W$ and a simplicial isomorphism $i: f({\cal A}) \ra f'({\cal A})$ such that $f' = i \circ f \circ w$. Note that this is a coarser equivalence relation than the equivalence relation obtained by the action of $(G \times W)$.

\bthm \label{thm:limites_a2} If ${\cal I}$ has type $A_2$, then the space $\Mor({\cal A},{\cal I})$ has $7$ equivalence classes, whose representatives are
\beq &\id : {\cal A} \ra {\cal A}, \hspace{0.5cm} p_1 : {\cal A} \ra \alpha_1 , \hspace{0.5cm} p_2 \circ p_1 : {\cal A} \ra \alpha_1 \cap \alpha_2 = \{C_2,C_3\},&\\
&p_6 \circ p_1  : {\cal A} \ra \alpha_1 \cap \alpha_6 = \{C_1,C_2\}, \hspace{0.5cm} p_3 \circ p_2 \circ p_1 : {\cal A} \ra \{C_3\}, \hspace{0.5cm} t_{12} \mbox{ and } t_{23}.&\eeq
According to Lemma~\ref{lem:proj}, we deduce that $\ov{\Mor_{inj}({\cal A},{\cal I})}^g = \Mor({\cal A},{\cal I})$.
\ethm

\bp Let $f \in \Mor({\cal A},{\cal I}) \bs \Mor_{inj}({\cal A},{\cal I})$.
\bit
\item If $f({\cal A})$ is included in an apartment, since ${\cal I}$ is topologically transitive, up to postcomposing with an element of $G$, we can assume that $f({\cal A}) \subset {\cal A}$. Since $f$ is not injective, assume that $f({\cal A}) \subset \alpha_1$.
\bit
\item If $f({\cal A}) = \alpha_1$, up to precomposing with an element of $W$, we can assume that $f=p_1$. Applying Lemma~\ref{lem:proj} to $\id : {\cal A} \hookrightarrow {\cal I}$, we get $f \in \ov{\Mor_{inj}({\cal A},{\cal I})}^g$.
\item If $f({\cal A})$ consists of two edges, then $f({\cal A})=\{C_2,C_3\}$ or $f({\cal A})=\{C_1,C_2\}$. For instance, assume $f({\cal A})=\{C_2,C_3\}$ : up to precomposing with an element of $W$, we can assume that $\Card f^{-1}(C_2) = 4$, and that $f = p_2 \circ p_1$. Applying Lemma~\ref{lem:proj} to $p_1$, we get $f \in \ov{\Mor_{inj}({\cal A},{\cal I})}^g$.
\item If $f({\cal A})$ is just one edge, assume $f({\cal A})=\{C_3\}$. Then $f=p_3 \circ p_2 \circ p_1$, and applying Lemma~\ref{lem:proj} to $p_2 \circ p_1$, we get $f \in \ov{\Mor_{inj}({\cal A},{\cal I})}^g$.
\eit
\item If $f({\cal A})$ is not included in an apartment, then $f$ is a marked tripod. Assume for instance that $f$ has the same type as $t_{12}$, and assume up to postcomposing with an element of $G$ that $f=t_{12}$ : that is $f(C_6)=f(C_1)=C_1$, $f(C_2)=f(C_3)=C_2$ and $f(C_4)=f(C_5)=C$, where $C$ is an edge in ${\cal I}$ containing $x_{12}$, different from $C_1$ and $C_2$.

Since ${\cal I}$ is topologically Moufang and has no isolated edge, we can find $(g_n)_{n \in \NNN}$ be a sequence in the root group $U_1 \subset G$ such that the sequence of edges $(g_n \cdot C_6)_{n \in \NNN}$ converges to $C_1$. Then the sequence of edges $(g_n \cdot C_5)_{n \in \NNN}$ converges to $C_2$. So the sequence of vertices $(g_n \cdot x_{56})_{n \in \NNN}$ converges to $x_{12}$.

Since ${\cal I}$ is topologically strongly transitive, the compact topology on the space of vertices of ${\cal I}$ with the same type as $x_{12}$ is the same as the topology induced by the topology of $G$-homogeneous space. Hence there exists a sequence $(g'_n)_{n \in \NNN}$ in $G$, converging to $e$, such that for every $n \in \NNN$ we have $g_n \cdot x_{56} = g'_n \cdot x_{12}$. Hence the sequence of stars $(\st(g_n \cdot x_{56}) = \st(g'_n \cdot x_{12}))_{n \in \NNN}$ converges to $\st(x_{12})$ in the space of closed subsets of ${\cal I}^{\<1\>}$.

In particular, for every $n \in \NNN$, there exists an edge $C'_n \in \st(g_n \cdot x_{56})$ such that the sequence of edges $(C'_n)_{n \in \NNN}$ converges to $C \in \st(x_{12})$. Since ${\cal I}$ is topologically Moufang, for every $n \in \NNN$, we can find $h_n \in U_6$ be such that $g_nh_n \cdot C_5 = C'_n$. It follows that the sequence of marked apartments $(g_nh_n : {\cal A} \ra {\cal I})_{n \in \NNN}$ converges to $f$ in $\Mor({\cal A},{\cal I})$, hence $f \in \ov{\Mor_{inj}({\cal A},{\cal I})}^g$. 
\eit
\ep

Applying Theorem~\ref{thm:compactification_plats}, we obtain the following result.

\bcor Let $X$ be the symmetric space of non-compact type or the Bruhat-Tits building of the group $\PGL(3)$ over a local field. Then the visual limits of divergent sequences of maximal flats in $X$ are the connected unions of $1$, $2$ or $3$ Weyl chambers in an apartment, and the tripods of the two types. \ecor

Note that in the case where ${\cal I}$ is the spherical topological building associated with the group $\PGL(3)$ over a local field, then each equivalence class in Theorem~\ref{thm:limites_a2} is an orbit of the $(G \times W)$-action.

\subsubsection{Type $C_2$ : topologically Moufang quadrangles}

\begin{figure}[!h]
\begin{center}
\includegraphics[width=15cm]{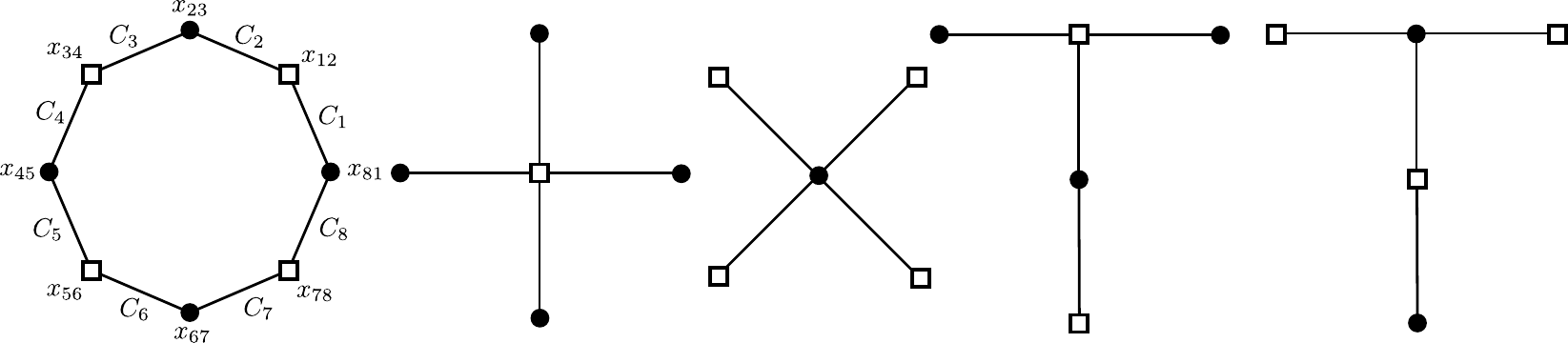}
\caption{The type $C_2$ apartment, the quadripods and the T-shapes}
\label{fig:c2}
\end{center}
\end{figure}

Assume here that $p=4$. Call (marked) \df{quadripod} of ${\cal I}$ any morphism $f : {\cal A} \ra {\cal I}$ whose image consists of four pairwise distinct and pairwise intersecting edges. Call \df{type} of a quadripod the type of this common intersection (see Figure~\ref{fig:c2}).

\medskip

For example, let $C,C'$ be two edges in ${\cal I}$ containing $x_{12}$, both different from $C_1$ and $C_2$. Then $f : {\cal A} \ra {\cal I}$, defined by $f(C_8)=f(C_1)=C_1$, $f(C_2)=f(C_3)=C_2$, $f(C_4)=f(C_5)=C$ and $f(C_6)=f(C_7)=C'$, is a quadripod.

\bigskip

Call (marked) \df{T-shape} of ${\cal I}$ any morphism $f : {\cal A} \ra {\cal I}$ whose image consists of four pairwise different edges $\{D_1,D'_1,D_2,D_3\}$, which are pairwise adjacent except that $D_3$ is not adjacent to $D_1$ nor $D'_1$ (see Figure~\ref{fig:c2}). Call \df{type} of a T-shape the type of the intersection of $D_1,D'_1$ and $D_2$.

\medskip

For example, fix $C'_1 \in \st(x_{12}) \bs \{C_1,C_2\}$, and consider $t_{12} : {\cal A} \ra {\cal I}$ defined by $t_{12}(C_8)=t_{12}(C_1)=C_1$, $t_{12}(C_2)=t_{12}(C_5)=C_2$, $t_{12}(C_3)=t_{12}(C_4)=C_3$ and $t_{12}(C_6)=t_{12}(C_7)=C'_1$, it
is a T-shape of one type. Similarly, fix $C'_2 \in \st(x_{23}) \bs \{C_2,C_3\}$, and consider $t_{23} : {\cal A} \ra {\cal I}$ defined by $t_{23}(C_1)=t_{23}(C_2)=C_2$, $t_{23}(C_3)=t_{23}(C_6)=C_3$, $t_{23}(C_4)=t_{23}(C_5)=C_4$ and $t_{23}(C_7)=t_{23}(C_8)=C'_2$, it is a T-shape of the other type.

Consider also degenerate quadripods and T-shapes, whose support are tripods. Define $t'_{12} : {\cal A} \ra \{C_1,C'_1,C_2\}$ and $t''_{12} : {\cal A} \ra \{C_1,C'_1,C_2\}$ by
\beq &t'_{12}(C_8)=t'_{12}(C_1)=C_1, \hspace{0.5cm} t'_{12}(C_2)=t'_{12}(C_3)=C_2,&\\
 &t'_{12}(C_4)=t'_{12}(C_5)=t'_{12}(C_6)=t'_{12}(C_7)=C'_1,&\\
&t''_{12}(C_8)=t''_{12}(C_1)=C_1, \hspace{0.5cm} t''_{12}(C_4)=t''_{12}(C_5)=C_2,&\\
&t''_{12}(C_2)=t''_{12}(C_3)=t''_{12}(C_6)=t''_{12}(C_7)=C'_1.&\eeq
Similarly, define $t'_{23} : {\cal A} \ra \{C_2,C'_2,C_3\}$ and $t''_{23} : {\cal A} \ra \{C_2,C'_2,C_3\}$ by
\beq &t'_{23}(C_1)=t'_{23}(C_2)=C_2, \hspace{0.5cm} t'_{23}(C_3)=t'_{23}(C_4)=C_3,&\\
&t'_{23}(C_5)=t'_{23}(C_6)=t'_{23}(C_7)=t'_{23}(C_8)=C'_2,&\\
&t''_{23}(C_1)=t''_{23}(C_2)=C_2, \hspace{0.5cm} t''_{23}(C_5)=t''_{23}(C_6)=C_3,&\\
&t''_{23}(C_3)=t''_{23}(C_4)=t''_{23}(C_7)=t''_{23}(C_8)=C'_2.&\eeq

Finally, define the morphism $f_0 : {\cal A} \ra \{C_1,C_2,C_3\}$ by
$$ f_0(C_8)=f_0(C_1)=C_1, f_0(C_2)=f_0(C_3)=f_0(C_4)=f_0(C_7)=C_2, f_0(C_5)=f_0(C_6)=C_3.$$

Recall that two elements $f$, $f'$ in $\Mor({\cal A},{\cal I})$ are equivalent if there exist $w \in W$ and a simplicial isomorphism $i: f({\cal A}) \ra f'({\cal A})$ such that $f' = i \circ f \circ w$.

\bthm \label{thm:limites_c2} If ${\cal I}$ has type $C_2$, then the space $\Mor({\cal A},{\cal I}) \bs \{\mbox{quadripods}\}$ has $19$ equivalence classes, whose representatives are
\bit
\item $\id : {\cal A} \ra {\cal A}$,
\item $p_1 : {\cal A} \ra \alpha_1$, $p_2 : {\cal A} \ra \alpha_2$,
\item $p_2 \circ p_1 : {\cal A} \ra \alpha_1 \cap \alpha_2 = \{C_2,C_3,C_4\}$, $p_8 \circ p_1 : {\cal A} \ra \alpha_1 \cap \alpha_8 = \{C_1,C_2,C_3\}$,
\item $f_0$,
\item $p_3 \circ p_2 \circ p_1 : {\cal A} \ra \{C_3,C_4\}$, $p_4 \circ p_3 \circ p_2 : {\cal A} \ra \{C_4,C_5\}$,
\item $p_3 \circ p_1 : {\cal A} \ra \{C_3,C_4\}$, $p_4 \circ p_2 : {\cal A} \ra \{C_4,C_5\}$,
\item $p_8 \circ p_2 \circ p_1 : {\cal A} \ra \{C_2,C_3\}$, $p_1 \circ p_3 \circ p_2 : {\cal A} \ra \{C_3,C_4\}$,
\item $p_4 \circ p_3 \circ p_2 \circ p_1 : {\cal A} \ra \{C_4\}$,
\item the T-shapes $t_{12}$ and $t_{23}$,
\item $t'_{12}$, $t'_{23}$,
\item $t''_{12}$ and $t''_{23}$.
\eit
The space $\ov{\Mor_{inj}({\cal A},{\cal I})}^g$ contains all these classes, except possibly the two classes of $t''_{12}$ and $t''_{23}$.
\ethm

\bp Let $f \in \Mor({\cal A},{\cal I}) \bs \Mor_{inj}({\cal A},{\cal I})$, which is not a quadripod. Recall that the Tits metric on ${\cal I}$ is here the length metric such that each edge has length $\frac{\pi}{4}$.
\bit

\item If $f({\cal A})$ is not included in an apartment, then the Tits diameter of $f({\cal A})$ is equal to $\frac{\pi}{2}$ or $\frac{3\pi}{4}$. Assume first that the diameter is $\frac{3\pi}{4}$, then $f$ is a T-shape. Up to the action of $(G \times W)$, we may assume that $f=t_{12}$ or $f=t_{23}$. For instance, assume that $f=t_{12}$.

Since ${\cal I}$ is topologically Moufang and has no isolated edge, let $(g_n)_{n \in \NNN}$ be a sequence in the root group $U_1 \subset G$ such that the sequence of edges $(g_n \cdot C_8)_{n \in \NNN}$ converges to $C_1$. Then the sequence of edges $(g_n \cdot C_7)_{n \in \NNN}$ converges to $C_2$. So the sequence of vertices $(g_n \cdot x_{78})_{n \in \NNN}$ converges to $x_{12}$.

Since ${\cal I}$ is topologically strongly transitive, the sequence of stars $(\st(g_n \cdot x_{78}))_{n \in \NNN}$ converges to $\st(x_{12})$ in the space of closed subsets of ${\cal I}^{\<1\>}$. In particular, for every $n \in \NNN$, there exists an edge $D_n \in \st(g_n \cdot x_{78})$ such that the sequence of edges $(D_n)_{n \in \NNN}$ converges to $C'_1 \in \st(x_{12})$. Since ${\cal I}$ is topologically Moufang, for every $n \in \NNN$, let $h_n \in U_8$ be such that $g_nh_n \cdot C_7 = D_n$. It follows that the sequence of marked apartments $(g_nh_n : {\cal A} \ra {\cal I})_{n \in \NNN}$ converges to $t_{12}$ in $\Mor({\cal A},{\cal I})$, hence $t_{12} \in \ov{\Mor_{inj}({\cal A},{\cal I})}^g$.

\item If the Tits diameter of $f({\cal A})$ is $\frac{\pi}{2}$, since $f$ is assumed not to be a quadripod, then $f({\cal A})$ is a tripod. Up to postcomposing with an element of $G$, assume that $f({\cal A})=\{C_1,C_2,C'_1\}$ or $f({\cal A})=\{C_2,C'_2,C_3\}$: assume the former. Then up to precomposing with an element of $W$, we may assume that $f=t'_{12}$ or $f=t''_{12}$. The theorem does not say anything about $t''_{12}$, so assume that $f=t'_{12}$.

Let $(D_n)_{n \in \NNN}$ be a sequence of edges in $\st(x_{23}) \bs \{C_2\}$ converging to $C_2$. For all $n \in \NNN$, let $t_n$ be the $T$-shape with image $\{C_1,C_2,D_n,C'_1\}$. The sequence $(t_n)_{n \in \NNN}$ converges to $t'_{12}$ in $\Mor_{inj}({\cal A},{\cal I})$, and $t_n \in \ov{\Mor_{inj}({\cal A},{\cal I})}^g$ by the previous case, so $t'_{12} \in \ov{\Mor_{inj}({\cal A},{\cal I})}^g$.

\item If $f({\cal A})$ is included in an apartment, up to postcomposing with an element of $G$, assume that $f({\cal A}) \subset {\cal A}$. Since $f$ is not injective, assume that $f({\cal A}) \subset \alpha_1$ or $f({\cal A}) \subset \alpha_2$.
\bit
\item If $f({\cal A}) = \alpha_1$ or $f({\cal A}) = \alpha_2$, up to precomposing with an element of $W$, we can assume that $f=p_1$ or $f=p_2$. According to lemma~\ref{lem:proj}, $f \in \ov{\Mor_{inj}({\cal A},{\cal I})}^g$.
\item If $f({\cal A})$ consists of three consecutive edges, we may assume that $f({\cal A})=\{C_1,C_2,C_3\}$. If $\Card f^{-1}(C_3) = 4$, then $f = p_8 \circ p_1$ and according to Lemma~\ref{lem:proj}, $f \in \ov{\Mor_{inj}({\cal A},{\cal I})}^g$. If $\Card f^{-1}(C_1) = 4$, then up to precomposing with an element $w$ of $W$, we can assume that $\Card (fw)^{-1}(C_1) = 4$. And if $\Card f^{-1}(C_2) = 4$, then $f=f_0$.

Let $(D_n)_{n \in \NNN}$ be a sequence of edges in $\st(x_{12}) \bs \{C_1,C_2\}$ converging to $C_2$. For all $n \in \NNN$, let $t_n$ be the $T$-shape with image $\{C_1,C_2,C_3,D_n\}$. The sequence $(t_n)_{n \in \NNN}$ converges to $f_0$ in $\Mor_{inj}({\cal A},{\cal I})$, and $t_n \in \ov{\Mor_{inj}({\cal A},{\cal I})}^g$ so $f_0 \in \ov{\Mor_{inj}({\cal A},{\cal I})}^g$.

\item If $f({\cal A})$ consists of two adjacent edges, we may assume that $f({\cal A})=\{C_3,C_4\}$ or $f({\cal A})=\{C_2,C_3\}$. For instance, assume $f({\cal A})=\{C_3,C_4\}$. If $\Card f^{-1}(C_3) = \Card f^{-1}(C_4) = 4$, then $f = p_3 \circ p_1$ or $f=p_1 \circ p_3 \circ p_2$, so according to Lemma~\ref{lem:proj}, $f \in \ov{\Mor_{inj}({\cal A},{\cal I})}^g$. If $\Card f^{-1}(C_3) = 6$ then $f=p_3 \circ p_2 \circ p_1$, so according to Lemma~\ref{lem:proj}, $f \in \ov{\Mor_{inj}({\cal A},{\cal I})}^g$.
\item If $f({\cal A})$ is just one edge, assume $f({\cal A})=\{C_4\}$. Then $f=p_4 \circ p_3 \circ p_2 \circ p_1$, and according to Lemma~\ref{lem:proj}, $f \in \ov{\Mor_{inj}({\cal A},{\cal I})}^g$. 

\eit
\eit
\ep

\subsubsection{Type $C_2$ : the classical case}

We will describe completely $\ov{\Mor_{inj}({\cal A},{\cal I})}^g$ in the classical case (see~\cite{tits_spherical_type} and \cite{moufang}).

\bigskip

Let $\KKK$ be a local field of characteristic different from $2$, or a quaternion algebra over a local field of characteristic different from $2$, and denote its center by $\ZZZ(\KKK)$. Denote by $\sigma : x \mapsto \ov{x}$ an involutive $\ZZZ(\KKK)$-antiautomorphism of $\KKK$, possibly trivial. Denote by $\KKK^\sigma = \{x \in \KKK \,:\, \ov{x}=x\}$ the $\ZZZ(\KKK)$-vector subspace of fixed points of $\sigma$ in $\KKK$.

If $\KKK$ is a quaternion algebra different from the Hamilton quaternion division algebra $\HHH$, we will assume that $\sigma$ is equal to $\id_\KKK$ or to the standard quaternionic involution $\sigma_0$. If we denote by $(1,i,j,k=ij)$ a standard $\ZZZ(\KKK)$-basis of the quaternion algebra $\KKK$, then $\sigma_0$ is defined by
$$ \sigma_0(i)=-i, \sigma_0(j)=-j \mbox{ and } \sigma_0(k)=-k.$$

If $\KKK=\HHH$ is the Hamilton quaternion division algebra, then up to automorphism there is only one involution different from the standard one: we will denote by $\tau$ such an involution, defined by
$$ \tau(i)=i, \tau(j)=-j \mbox{ and } \tau(k)=k.$$
So if $\KKK=\HHH$ is the Hamilton quaternion division algebra, we will assume (without loss of generality) that $\sigma$ is equal to $\id_\KKK$, to the standard quaternionic involution $\sigma_0$ or to $\tau$.

Let $V$ be a finite-dimensional right vector space over $\KKK$ of dimension at least $5$, and let $q$ be a nondegenerate $\sigma$-Hermitian form on $V$. It means that the associated form $\varphi : V \times V \ra \KKK$ is $\sigma$-sesquilinear, and $\sigma$-Hermitian symmetric (see~\cite{bourbaki_algebre9} and \cite{tits_spherical_type}):
$$ \forall v,w \in V, \forall x,y \in \KKK,\, \varphi(vx,wy)=\sigma(x)\varphi(v,w)y \mbox{ and } \varphi(w,v) = \sigma\left(\varphi(v,w)\right).$$
Assume further that the Witt index of $q$ is $2$, which means that all the maximal totally isotropic subspaces of $V$ have the same dimension $2$.

Remark that if $\sigma=\id_\KKK$, the existence of $q$ implies that $\KKK$ is commutative.

\bigskip

Let ${\cal I}$ be the flag complex of totally isotropic subspaces of $V$: it is a spherical building of type $C_2$, called \df{classical}. It is always thick, except when $\sigma=\id_\KKK$ and $\dim V=4$.

\bigskip

The projective space $\PPP(V)$ and the Grassmannian of $2$-planes $\Gr_2(V)$ are naturally endowed with compact, non-discrete topologies. Consider the topology on ${\cal I}^{\<1\>}$ induced by the product topology on $\PPP(V) \times \Gr_2(V)$: it turns ${\cal I}$ into a compact non-discrete topological spherical building. Let $G=\PU(q)$ be the projective isotropy group of the Hermitian form $q$, it acts strongly transitively on ${\cal I}$.

In the Archimedean case, this setting includes the following real noncompact simple Lie groups of real rank $2$ and of type $C_2$:
\bit
\item $\PO(2,n)$ for $n \geq 3$, when $\KKK=\RRR$ and $\dim V = 2+n$.
\item $\PO(5,\CCC)$, when $\KKK=\CCC$, $\sigma=\id_\CCC$ and $\dim V = 5$.
\item $\PU(2,n)$ for $n \geq 3$, when $\KKK=\CCC$, $\sigma \neq \id_\CCC$ and $\dim V = 2+n$.
\item $\PSp(2,n)$ for $n \geq 3$, when $\KKK=\HHH$, $\sigma = \sigma_0$ and $\dim V = 2+n$.
\item $\PO^*(10)$, when $\KKK=\HHH$, $\sigma = \tau$ and $\dim V = 5$.
\eit

\bigskip

Let us now describe the standard apartment ${\cal A}$ of ${\cal I}$. If $E$ is a subset of $V$, denote by $\<E\>$ the right $\KKK$-vector subspace spanned by $E$. Fix $(e_1,e_2,e_3,e_4)$ a free family of isotropic vectors of $V$ such that $\<e_1,e_3\>$ and $\<e_2,e_4\>$ are orthogonal hyperbolic planes, normalised in such a way that
$$ \varphi(e_1,e_3) = \varphi(e_2,e_4) = 1.$$
Let ${\cal A}$ denote the apartment of ${\cal I}$ whose \df{line-type} vertices are the four isotropic lines $x_{81}=\<e_1\>$, $x_{23}=\<e_2\>$, $x_{45}=\<e_3\>$ and $x_{67}=\<e_4\>$, and whose \df{plane-type} vertices are the four isotropic planes $x_{12}=\<e_1,e_2\>$, $x_{34}=\<e_2,e_3\>$, $x_{56}=\<e_3,e_4\>$ and $x_{78}=\<e_4,e_1\>$.

\bigskip

Let $q_{an}$ denote the restriction of $q$ to the orthogonal complement $V'$ of $\<e_1,e_2,e_3,e_4\>$ in $V$: since the Witt index of $q$ is $2$, the Hermitian form $q_{an}$ is anisotropic.

\bigskip

If $f$ is a line-type quadripod, call \df{dimension} of $f$ the dimension of the right $\KKK$-vector space spanned by the four plane-type vertices. Let $p_1,p_2,p_3,p_4$ denote the four cyclically-ordered isotropic $2$-planes of $f$, with $p_1 \cap p_2 \cap p_3 \cap p_4 = \ell$, an isotropic line. Say $f$ is \df{symmetric} if its dimension is at most $4$ and there exists a $\sigma$-semilinear semiisometric involution $s$ of $V$, such that $s(\ell)=\ell$ and $\forall i \in \lb 1,4 \rb$, $s(p_i)=p_{i+2}$.

Note that $f$ is symmetric if and only if for all $i \in \lb 1,4\rb$, there exists $v_i \in p_i \bs \ell$ such that $v_1+v_2+v_3+v_4=0$ and $\varphi(v_2,v_3) = \varphi(v_1,v_4)$, or equivalently $\varphi(v_1,v_2) = \varphi(v_4,v_3)$.

The notion of symmetric line-type quadripod is mostly important in the Hermitian case $\sigma \neq \id_\KKK$, according to the following result.

\blem \label{lem:orthogonal_symmetric} If $\sigma = \id_\KKK$ and $f$ is a line-type quadripod of dimension at most $4$, then $f$ is symmetric. \elem

\bp
Since $f$ is a (nondegenerate) quadripod, its dimension is at least $4$, so let us assume $f$ has dimension $4$. According to Witt's theorem and up to postcomposing with $G$, we may assume that the central line-type vertex of $f$ is $\<e_1\>$, that $f(\<e_1,e_2\>)=\<e_1,e_2\>$ and that $f(\<e_1,e_4\>)=\<e_1,e_4\>$. The $\KKK$-vector subspace $E$ of $V$ spanned by the plane-type vertices of $f$ has dimension $4$, contains $\<e_1,e_2,e_4\>$ and is included in the orthogonal complement of $\<e_1\>$, hence it is equal to $E=\<e_1,e_2,e_4,e_5\>$, where $e_5 \neq 0$ is orthogonal to $\<e_1,e_2,e_3,e_4\>$. Let $r \in \KKK \bs \{0\}$ such that $\varphi(e_5,e_5)=-2r$. Up to multiplying $e_2$, $e_4$ and $e_5$ by scalars, we may further assume that $f(\<e_2,e_3\>)=\<e_1,e_2+e_4 r+e_5\>$. And there exists $a \in \KKK$ such that $f(\<e_3,e_4\>)=\<e_1,e_2+e_4 a^2r+e_5 a\>$. The vectors $v_1=e_2(1-a)$, $v_2=(e_2+e_4 r+e_5)a$, $v_3=(e_2+e_4 r+e_5)(-1)$ and $v_4=e_4r(a^2-a)$ then satisfy $v_1+v_2+v_3+v_4=0$ and $\varphi(v_2,v_3) = \varphi(v_1,v_4)=-ar(a-1)^2$, so $f$ is symmetric.
\ep

If $f$ is a plane-type quadripod and $\KKK$ is commutative, call \df{cross-ratio} of $f$ the cross-ratio of the
four line-type vertices of $f$, which are all included in the central plane-type vertex of $f$.

\bthm \label{thm:limites_c2_classical} Suppose that ${\cal I}$ is of type $C_2$ and classical, and fix a morphism $f : {\cal A} \ra {\cal I}$.
\bit
\item If $f$ is not a quadripod, then $f \in \ov{\Mor_{inj}({\cal A},{\cal I})}^g$.
\item If $f$ is a line-type quadripod, then $f \in \ov{\Mor_{inj}({\cal A},{\cal I})}^g$ if and only if $f$ is symmetric.
\item If $f$ is a plane-type quadripod, then $f \in \ov{\Mor_{inj}({\cal A},{\cal I})}^g$ if and only if $\sigma \neq \id_\KKK$, or $\sigma=\id_\KKK$ and $1-c \in \{q_{an}(v)q_{an}(v')^{-1} \,:\, v,v' \in V' \bs \{0\}\}$, where $c$ is the cross-ratio of $f$.
\eit
Assume $\sigma=\id_\KKK$. If $\dim V = 5$, then all line-type quadripods belong to $\ov{\Mor_{inj}({\cal A},{\cal I})}^g$, and a plane-type quadripod with cross-ratio $c$ belongs to $\ov{\Mor_{inj}({\cal A},{\cal I})}^g$ if and only if $1-c$ is a square in $\KKK$. So if $\KKK=\RRR$, then the condition on $c$ reduces to $c \leq 1$, and if further $\KKK=\CCC$ then the condition on $c$ is void. In particular if $\KKK=\CCC$, we have
$$\Mor({\cal A},{\cal I}) = \ov{\Mor_{inj}({\cal A},{\cal I})}^g.$$
\ethm

\bp \bit
\item Consider the orthogonal case $\sigma=\id_\KKK$. Then $\KKK$ is commutative, hence $\KKK$ is a local field.

Let us first focus on quadripods.
\bit
\item Let $f$ be a line-type quadripod. If ${\cal A}'$ is an apartment in ${\cal I}$, the $\KKK$-vector subspace spanned by the four plane-type vertices and the four line-type vertices of ${\cal A}'$ has dimension $4$. Hence if $f \in \ov{\Mor_{inj}({\cal A},{\cal I})}^g$, the dimension of $f$ is at most $4$ (by lower semicontinuity of the dimension). Conversely, let us assume that the dimension of $f$ is at most $4$: since $f$ is a (nondegenerate) quadripod, its dimension is at least $4$, so let us assume it is equal to $4$. Then, according to Witt's theorem and up to postcomposing with $G$, we may assume that the central line-type vertex of $f$ is $\<e_1\>$, that $f(\<e_1,e_2\>)=\<e_1,e_2\>$ and that $f(\<e_1,e_4\>)=\<e_1,e_4\>$. The $\KKK$-vector subspace $E$ of $V$ spanned by the plane-type vertices of $f$ has dimension $4$, contains $\<e_1,e_2,e_4\>$ and is included in the orthogonal complement of $\<e_1\>$, hence it is equal to
$E=\<e_1,e_2,e_4,e_5\>$, where $e_5$ is orthogonal to $\<e_1,e_2,e_3,e_4\>$. Let $r \in \KKK \bs \{0\}$ such that $\varphi(e_5,e_5)=-2r$. Up to multiplying $e_2$, $e_4$ and $e_5$ by scalars, we may further assume that $f(\<e_2,e_3\>)=\<e_1,e_2+e_4 r+e_5\>$. And there exists $a \in \KKK \bs \{0,1\}$ such that $f(\<e_3,e_4\>)=\<e_1,e_2+e_4 a^2r+e_5 a\>$.

\bigskip
Let $(x_n)_{n \in \NNN }$ be a sequence in $\KKK$ going to infinity. Let us define the marked apartment $f_n \in \Mor_{inj}({\cal A},{\cal I})$ by
$f_n(\<e_i\>) = \<\widetilde{f_n}(e_i)\>$ for all $i \in \lb 1,4 \rb$, where
\beq &\widetilde{f_n}(e_1)=e_1, \hspace{0.5cm}  \widetilde{f_n}(e_2) = -e_1arx_n+e_2, \hspace{0.5cm}  \widetilde{f_n}(e_4) = -e_1 x_n+e_4& \\
& \mbox{ and } \widetilde{f_n}(e_3) = e_1 ar(a-1)x_n^2 + e_2 x_n + e_3 + e_4 arx_n + e_5ax_n.&\eeq
Each of these vectors is isotropic :
\beq &q(\widetilde{f_n}(e_1))=0, \hspace{0.5cm} q(\widetilde{f_n}(e_2))=0, \hspace{0.5cm} q(\widetilde{f_n}(e_4))=0& \\
&\mbox{ and } q(\widetilde{f_n}(e_3)) = 2ar(a-1)x_n^2 + 2arx_n^2 + (-2r)a^2x_n^2=0.&\eeq
Furthermore $\varphi(\widetilde{f_n}(e_1),\widetilde{f_n}(e_2)) = \varphi(\widetilde{f_n}(e_4),\widetilde{f_n}(e_1)) = 0$, and
\beq \varphi(\widetilde{f_n}(e_2),\widetilde{f_n}(e_3)) &=& -arx_n +arx_n= 0 \\
 \varphi(\widetilde{f_n}(e_3),\widetilde{f_n}(e_4)) &=& x_n-x_n= 0.\eeq
Hence $f_n$ is a marked apartment of ${\cal I}$. Let us show that $(f_n)_{n \in \NNN}$ converges to $f$.

First of all, for all $i \in \lb 1,4\rb$, the sequence $(f_n(\<e_i\>))_{n \in \NNN}$ converges to $\<e_1\> = f(\<e_1\>)$. For all $n \in \NNN$, we have $f_n(\<e_1,e_2\>) = \<e_1,e_2\> =f(\<e_1,e_2\>)$ and $f_n(\<e_1,e_4\>) = \<e_1,e_4\> =f(\<e_1,e_4\>)$.

Every accumulation point of the sequence of planes $(f_n(\<e_2,e_3\>))_{n \in \NNN}$ contains the line
\beq \liml_{n \ra +\infty} \<\widetilde{f_n}(e_3)+\widetilde{f_n}(e_2)(arx_n)^{-1}ar(a-1)x_n^2\> &=&\\
\liml_{n \ra +\infty} \<e_2(x_n + (a-1)x_n)+e_3+ e_4 arx_n + e_5 ax_n\> &=&\\
\<e_2+e_4 r+ e_5\> \subset f(\<e_2,e_3\>).&&\eeq
Hence the sequence of planes $(f_n(\<e_2,e_3\>))_{n \in \NNN}$ converges to $f(\<e_2,e_3\>)$.

Similarly, every accumulation point of the sequence of planes $(f_n(\<e_3,e_4\>))_{n \in \NNN}$ contains the line
\beq \liml_{n \ra +\infty} \<f_n(e_3)+f_n(e_4)x_n^{-1}ar(a-1)x_n^2\> &=&\\
\liml_{n \ra +\infty} \<e_2x_n+e_3+ e_4(arx_n + ar(a-1)x_n) + e_5ax_n\> &=&\\
\<e_2+e_4 a^2r + e_5a)\> \subset f(\<e_3,e_4\>).&&\eeq
Hence the sequence of planes $(f_n(\<e_3,e_4\>))_{n \in \NNN}$ converges to $f(\<e_3,e_4\>)$.

To conclude this case, the sequence of marked apartments $(f_n)_{n \in \NNN}$ converges to $f$, so $f \in \ov{\Mor_{inj}({\cal A},{\cal I})}^g$.

\item Let $f$ be a plane-type quadripod. Then, according to Witt's theorem and up to postcomposing with $G$, we may assume that the central plane-type vertex of $f$ is $\<e_1,e_2\>$, that $f(\<e_1\>)=\<e_1\>$ and that $f(\<e_2\>)=\<e_2\>$. Up to multiplying $e_1$ and $e_2$ by scalars, we may further assume that if $c \in \KKK \bs \{0,1\}$ is the cross-ratio of $f$, then
$$ f(\<e_1\>)= \<e_1\>, \hspace{0.5cm} f(\<e_2\>)= \<e_2\>, \hspace{0.5cm} f(\<e_3\>)= \<e_1+e_2\> \mbox{ and } f(\<e_4\>)= \<e_1 c+e_2\>.$$

\smallskip

Assume first that there exist $v,v' \in V' \bs \{0\}$ such that $1-c=q(v)q(v')^{-1}$. Let $(x_n)_{n \in \NNN }$ be a sequence in $\KKK$ leaving every compact subset. Let us define the marked apartment $f_n \in \Mor_{inj}({\cal A},{\cal I})$ by $f_n(\<e_i\>) = \<\widetilde{f_n}(e_i)\>$ for all $i \in \lb 1,4 \rb$, where
\beq &\widetilde{f_n}(e_1)=e_1, \hspace{0.5cm} \widetilde{f_n}(e_2) = e_2,&\\
 &\widetilde{f_n}(e_3) = (e_1+e_2)(-\frac{1}{2}q(v+v')x_n^2)+e_3+(v+v')x_n& \mbox{ and }\\
&  \widetilde{f_n}(e_4) = (e_1c+e_2)(-\frac{1}{2}q(v')x_n^2)+e_4+v'x_n.&\eeq

Each of these vectors is isotropic :
\beq &q(\widetilde{f_n}(e_1))=0, \hspace{0.5cm} q(\widetilde{f_n}(e_2))=0,& \\
&q(\widetilde{f_n}(e_3))=2(-\frac{1}{2}q(v+v')x_n^2)+q(v+v')x_n^2=0& \mbox{and }\\
& q(\widetilde{f_n}(e_4))=2(-\frac{1}{2}q(v')x_n^2)+q(v')x_n^2=0.&\eeq
Furthermore $\varphi(\widetilde{f_n}(e_1),\widetilde{f_n}(e_2)) = \varphi(\widetilde{f_n}(e_2),\widetilde{f_n}(e_3)) = \varphi(\widetilde{f_n}(e_4),\widetilde{f_n}(e_1)) = 0$, and
\beq \varphi(\widetilde{f_n}(e_3),\widetilde{f_n}(e_4)) &=& -\frac{1}{2}q(v+v')x_n^2 - \frac{c}{2}q(v')x_n^2 + \varphi(v+v',v')x_n^2\\
&=& -\frac{1}{2}(q(v+v') + cq(v') -2\varphi(v+v',v'))x_n^2\\
&=& -\frac{1}{2}(q(v+v') + q(v')-q(v) -2\varphi(v+v',v'))x_n^2 = 0.\eeq
Hence $f_n$ is a marked apartment of ${\cal I}$. And it is clear that $(f_n)_{n \in \NNN}$ converges to $f$, so $f \in \ov{\Mor_{inj}({\cal A},{\cal I})}^g$.

\smallskip

Conversely, assume that $f \in \ov{\Mor_{inj}({\cal A},{\cal I})}^g$: let $(f_n)_{n \in \NNN}$ be a sequence of marked apartments converging to $f$. Then, according to Witt's theorem and up to postcomposing with $G$, we may assume that, for all $n \in \NNN$, we have $f_n(\<e_1\>)= \<e_1\>$ and $f_n(\<e_2\>)= \<e_2\>$. For all $n \in \NNN$, since $f_n$ is a marked apartment we know there exists $x_n, y_n, z_n, t_n \in \KKK$ and $u_n,v_n \in V'$ such that
$$ f_n(\<e_3\>)= \<e_1+e_2x_n+e_3z_n+u_n\> \mbox{ and } f_n(\<e_4\>)= \<e_1c+e_2y_n+e_4t_n+v_n\>.$$
As the sequence $(f_n)_{n \in \NNN}$ converges to $f$, we deduce that $x_n,y_n \ral{n \ra +\infty} 1$, that $z_n,t_n \ral{n \ra +\infty} 0$ and that $u_n,v_n \ral{n \ra +\infty} 0$. Furthermore $f_n(\<e_3\>)$ and $f_n(\<e_4\>)$ are isotropic and orthogonal to each other, so we conclude that $2z_n+q(u_n)=0$, that $2y_nt_n+q(v_n)=0$ and that $x_nt_n+cz_n+\varphi(u_n,v_n)=0$.
Hence
\beq c=(-x_nt_n-\varphi(u_n,v_n))z_n^{-1} &=& (2^{-1}x_ny_n^{-1}q(v_n)-\varphi(u_n,v_n))\left(-\frac{1}{2}q(u_n)\right)^{-1}, \\
(1-c)q(u_n) &=& q(u_n) + x_ny_n^{-1}q(v_n) - 2\varphi(u_n,v_n).\eeq
Since $q_{an}$ is anisotropic, $\sqrt{|q_{an}|}$ is equivalent to a norm on $V'$, so we can assume that, up to passing to a subsequence, one of the sequences $(q(v_n)q(u_n)^{-1})_{n \in \NNN}$ or $(q(v_n)q(u_n-v_n)^{-1})_{n \in \NNN}$ is bounded. As a consequence, $(1-c)q(u_n) \underset{n \ra +\infty}{\sim} q(u_n) + q(v_n) - 2\varphi(u_n,v_n) = q(u_n-v_n)$ so
$$ 1-c = \liml_{n \ra +\infty} q(u_n-v_n)q(u_n)^{-1}.$$
Since the set $\{q_{an}(v)q_{an}(v')^{-1} \,:\, v,v' \in V' \bs \{0\}\}$ is stable under multiplication by squares in $\KKK \bs \{0\}$, it is closed in $\KKK \bs \{0\}$ since the subgroup of squares in the multiplicative group $\KKK \bs \{0\}$ is of finite index. Hence
$$1-c \in \{q_{an}(v)q_{an}(v')^{-1} \,:\, v,v' \in V' \bs \{0\}\}.$$
\eit

\bigskip

By Theorem~\ref{thm:limites_c2}, to prove Theorem~\ref{thm:limites_c2_classical} in the case $\sigma=\id_\KKK$ it remains to prove that the morphisms $t''_{12}$ and $t''_{23}$ belong to $\ov{\Mor_{inj}({\cal A},{\cal I})}^g$.

Firstly, the morphism $t''_{12}$ can be interpreted as a degenerate plane-type quadripod, whose cross-ratio is $0$. The set $\{q_{an}(v)q_{an}(v')^{-1} \,:\, v,v' \in V' \bs \{0\}\}$ is stable under multiplication by squares in $\KKK \bs \{0\}$, and since the subgroup of squares in the multiplicative group $\KKK \bs \{0\}$ is of finite index, then $1$ is the limit of a sequence $(c_n)_{n \in \NNN}$ in $\{q_{an}(v)q_{an}(v')^{-1} \,:\, v,v' \in V' \bs \{0\}\} \bs \{1\}$. Then there exists a sequence $(f_n)_{n \in \NNN}$ of plane-type quadripods, whose cross-ratios are $(c_n)_{n \in \NNN}$, which converges to $t''_{12}$. For all $n \in \NNN$, since $c_n \not\in \{0,1\}$ we know that $f_n$ is a nondegenerate plane-type quadripod, which then belongs to $\ov{\Mor_{inj}({\cal A},{\cal I})}^g$. Hence $t''_{12} \in \ov{\Mor_{inj}({\cal A},{\cal I})}^g$.

Secondly, the morphism $t''_{23}$ can be interpreted as a degenerate line-type quadripod, for which two opposite plane-type vertices coincide. Let $E$ denote the vector space spanned by the plane-type vertices of $t''_{23}$, it has dimension at most $4$. Then there exists a sequence $(f_n)_{n \in \NNN}$ of line-type quadripods, for which the vector space spanned by the plane-type vertices is $E$, which converges to $t''_{23}$. But we know that each $f_n$ has dimension equal to $\dim E \leq 4$ and belongs to $\ov{\Mor_{inj}({\cal A},{\cal I})}^g$, and so $t''_{23}$ belongs to $\ov{\Mor_{inj}({\cal A},{\cal I})}^g$.

\bigskip

Assume now that $\dim V=5$. Then since $q$ is nondegenerate, the dimension of any line-type quadripod is at most $4$. And $V'=\<e_1,e_2,e_3,e_4\>^\perp$ is a line, and so the set $\{q_{an}(v)q_{an}(v')^{-1} \,:\, v,v' \in V' \bs \{0\}\}$ is just the set of squares in $\KKK \bs \{0\}$.

\bigskip

\item Consider now the Hermitian case $\sigma \neq \id_\KKK$.

Let us first focus on quadripods.
\bit
\item Let $f$ be a line-type quadripod. Let us assume that the dimension of $f$ is at most $4$ and that $f$ is symmetric. Since $f$ is nondegenerate and symmetric, then its dimension is exactly $4$. According to Witt's theorem and up to postcomposing with $G$, we may assume that the central line-type vertex of $f$ is $\<e_1\>$, that $f(\<e_1,e_2\>)=\<e_1,e_2\>$ and that $f(\<e_1,e_4\>)=\<e_1,e_4\>$. The $\KKK$-vector subspace $E$ of $V$ spanned by the plane-type vertices of $f$ has dimension $4$, contains $\<e_1,e_2,e_4\>$ and is included in the orthogonal complement of $\<e_1\>$, hence it is equal to
$E=\<e_1,e_2,e_4,e_5\>$, where $e_5$ is orthogonal to $\<e_1,e_2,e_3,e_4\>$. Denote $\varphi(e_5,e_5)=r \in \KKK^\sigma \bs \{0\}$. There exist $a,b,c,d \in \KKK$ such that
$$ f(\<e_2,e_3\>)=\<e_1,e_2+e_4 a+e_5b\> \mbox{ and } f(\<e_3,e_4\>)=\<e_1,e_2+e_4 c+e_5 d\> ,$$
where $a+\ov{a}+\ov{b}rb=c+\ov{c}+\ov{d}rd=0$.

Since the dimension of $f$ is $4$, $b$ and $d$ are not both zero. And since $f$ is symmetric, neither $b$ nor $d$ can be zero, hence $a$ and $c$ are non-zero as well.

Let $(x_n)_{n \in \NNN }$ be a sequence in $\ZZZ(\KKK)^\sigma$ going to infinity. Let us define the marked apartment $f_n \in \Mor_{inj}({\cal A},{\cal I})$ by $f_n(\<e_i\>) = \<\widetilde{f_n}(e_i)\>$ for all $i \in \lb 1,4 \rb$, where
\beq &\widetilde{f_n}(e_1)=e_1, \hspace{0.5cm} \widetilde{f_n}(e_2) = e_1(-x_n)+e_2, \hspace{0.5cm} \widetilde{f_n}(e_4) = e_1 (-\ov{a^{-1}}\ov{b}\ov{d^{-1}}x_n)+e_4& \\
& \mbox{ and } \widetilde{f_n}(e_3) = e_1 (a^{-1}-d^{-1}ba^{-1})x_n^2 + e_2 d^{-1}ba^{-1}x_n + e_3 + e_4 x_n + e_5 ba^{-1}x_n.&\eeq
Each of these vectors is isotropic $q(\widetilde{f_n}(e_1))=q(\widetilde{f_n}(e_2))=q(\widetilde{f_n}(e_4))=0$ and
\beq q(\widetilde{f_n}(e_3)) &=& (\ov{a^{-1}}-\ov{a^{-1}}\ov{b}\ov{d^{-1}})x_n^2 + (a^{-1}-d^{-1}ba^{-1})x_n^2 \\
& &+\ \ov{a^{-1}}\ov{b}\ov{d^{-1}}x_n^2 + d^{-1}ba^{-1}x_n^2 + \ov{a^{-1}}\ov{b}rba^{-1}x_n^2 = 0.\eeq
Furthermore $\varphi(\widetilde{f_n}(e_1),\widetilde{f_n}(e_2)) = \varphi(\widetilde{f_n}(e_4),\widetilde{f_n}(e_1)) = 0$, and
\beq \varphi(\widetilde{f_n}(e_2),\widetilde{f_n}(e_3)) &=& -x_n +x_n= 0 \\
 \varphi(\widetilde{f_n}(e_3),\widetilde{f_n}(e_4)) &=& \ov{d^{-1}ba^{-1}}x_n + (-\ov{a^{-1}}\ov{b}\ov{d^{-1}}x_n) = 0.\eeq
Hence $f_n$ is a marked apartment of ${\cal I}$. Let us show that $(f_n)_{n \in \NNN}$ converges to $f$.

First of all, for all $i \in \lb 1,4\rb$, the sequence $(f_n(\<e_i\>))_{n \in \NNN}$ converges to $\<e_1\> = f(\<e_1\>)$. For all $n \in \NNN$, we have $f_n(\<e_1,e_2\>) = \<e_1,e_2\> =f(\<e_1,e_2\>)$ and $f_n(\<e_1,e_4\>) = \<e_1,e_4\> =f(\<e_1,e_4\>)$.

Every accumulation point of the sequence of planes $(f_n(\<e_2,e_3\>))_{n \in \NNN}$ contains the line
\beq \liml_{n \ra +\infty} \<\widetilde{f_n}(e_3)+\widetilde{f_n}(e_2)(a^{-1}-d^{-1}ba^{-1})x_n\> &=&\\
\liml_{n \ra +\infty} \<e_2(d^{-1}ba^{-1} + a^{-1}-d^{-1}ba^{-1})x_n+e_3+ e_4 x_n + e_5 ba^{-1}x_n\> &=&\\
\<e_2+e_4 a+ e_5 b\> \subset f(\<e_2,e_3\>).&&\eeq
Hence the sequence of planes $(f_n(\<e_2,e_3\>))_{n \in \NNN}$ converges to $f(\<e_2,e_3\>)$.

Similarly, every accumulation point of the sequence of planes $(f_n(\<e_3,e_4\>))_{n \in \NNN}$ contains the line
\beq \liml_{n \ra +\infty} \<\widetilde{f_n}(e_3)+\widetilde{f_n}(e_4)(\ov{a^{-1}}\ov{b}\ov{d^{-1}})^{-1}(a^{-1}-d^{-1}ba^{-1})x_n\> &=&\\
\liml_{n \ra +\infty} \<e_2 d^{-1}ba^{-1}x_n+e_3+ e_4(1+\ov{d}\ov{b^{-1}}\ov{a}(a^{-1}-d^{-1}ba^{-1}))x_n + e_5ba^{-1}x_n\> &=&\\
\<e_2+e_4(ab^{-1}d+\ov{d}\ov{b^{-1}}\ov{a}(b^{-1}d-1)) + e_5 d)\>.&&\eeq

So it remains to show that $ab^{-1}d+\ov{d}\ov{b^{-1}}\ov{a}(b^{-1}d-1)=c$. Let
\beq v_{12}&=&e_2(d^{-1}-b^{-1}) \in f(\<e_1,e_2\>) \bs \<e_1\>,\\
v_{23}& =& e_2 b^{-1}+e_4ab^{-1}+e_5 \in f(\<e_2,e_3\>) \bs \<e_1\>\\
v_{34} &=& -e_2 d^{-1}-e_4cd^{-1}-e_5 \in f(\<e_3,e_4\>) \bs \<e_1\>\\
\mbox{ and } v_{41}&=&e_4(cd^{-1}-ab^{-1}) \in f(\<e_4,e_1\>) \bs \<e_1\>,\eeq
they are such that $v_{12}+v_{23}+v_{34}+v_{41}=0$, so the symmetry assumption about $f$ tells us that
\beq \varphi(v_{12},v_{23}) &=& \varphi(v_{41},v_{34}) \\
(\ov{d}^{-1}-\ov{b}^{-1})ab^{-1} &=& (\ov{d}^{-1}\ov{c}-\ov{b}^{-1}\ov{a})(-d^{-1}) \\
\ov{d}^{-1}\ov{c}d^{-1} &=& \ov{b}^{-1}\ov{a}d^{-1} - \ov{d}^{-1}ab^{-1}+\ov{b}^{-1}ab^{-1} \\
c &=& ab^{-1}d + \ov{d}\ov{b}^{-1}\ov{a}b^{-1}d - \ov{d}\ov{b}^{-1}\ov{a}.\eeq
Hence the sequence $(f_n)_{n \in \NNN}$ converges to $f$.

\bigskip

Conversely, assume that $f \in \ov{\Mor_{inj}({\cal A},{\cal I})}^g$. If ${\cal A}'$ is an apartment in ${\cal I}$, the $\KKK$-vector space spanned by the four plane-type vertices of ${\cal A}'$ has dimension $4$. Hence the dimension of $f$ is at most $4$, and so we can choose $v_{12} \in f(\<e_1,e_2\>) \bs f(\<e_1\>)$, $v_{23} \in f(\<e_2,e_3\>) \bs f(\<e_1\>)$, $v_{34} \in f(\<e_3,e_4\>) \bs f(\<e_1\>)$ and $v_{41} \in f(\<e_4,e_1\>) \bs f(\<e_1\>)$ such that $v_{12}+v_{23}+v_{34}+v_{41}=0$.

Let $(f_n)_{n \in \NNN}$ be a sequence of marked apartments converging to $f$. For $i \in \ZZZ/4\ZZZ$, choose a sequence $(v^n_{i,i+1})_{n \in \NNN}$ in $(f_n(\<e_i,e_{i+1}\>))_{n \in \NNN}$ converging to $v_{i,i+1}$, such that for all $n \in \NNN$ we have $v^n_{12}+v^n_{23}+v^n_{34}+v^n_{41}=0$. Hence for all $n \in \NNN$, we have
$$ \varphi(v^n_{23},v^n_{43}) = \varphi(-v^n_{12},-v^n_{41}) = \varphi(v^n_{12},v^n_{41}).$$
Taking the limit $n \ra +\infty$, we get $\varphi(v_{23},v_{43}) = \varphi(v_{12},v_{41})$, so $f$ is symmetric.

\item Let $f$ be a plane-type quadripod. Then, according to Witt's theorem and up to postcomposing with $G$, we may assume that the central plane-type vertex of $f$ is $\<e_1,e_2\>$, that $f(\<e_1\>)=\<e_1\>$ and that $f(\<e_2\>)=\<e_2\>$. Up to multiplying $e_1$ and $e_2$ by scalars, we may further assume that there exists $c \in \KKK \bs \{0,1\}$ such that
$$ f(\<e_1\>)= \<e_1\>, f(\<e_2\>)= \<e_2\>, f(\<e_3\>)= \<e_1+e_2\> \mbox{ and } f(\<e_4\>)= \<e_1 c+e_2\>.$$
In fact, the scalar $c$ is the cross-ratio of $f$ (or an element of the cross-ratio if $\KKK$ is not commutative: see~\cite{baer} for instance). If $\KKK$ is a quaternion division algebra and $\sigma=\sigma_0$ is the standard quaternionic involution, then $c \ov{c} \in \KKK^\sigma=\ZZZ(\KKK)$, hence $\ov{c}$ commutes with $c$. If $\KKK=\HHH$ is the Hamilton quaternion division algebra and $\sigma=\tau$, then since $c$ is well-defined up to conjugation we may assume that $c \in \RRR \oplus \RRR j \subset \HHH$, hence $\ov{c}$ commutes with $c$.

\smallskip

Fix $e_5 \in \<e_1,e_2,e_3,e_4\>^\perp \bs \{0\}$, and denote $r=q(e_5) \in \KKK^\sigma \bs \{0\}$. If $\KKK^\sigma$ is not included in the center of $\KKK$, then by assumption $\KKK = \HHH$ and $\sigma = \tau$. Then $r$ is a square in $\HHH^\tau \bs \{0\}$, hence up to multiplying $e_5$ by a scalar we may assume that $r=1$, so that in any case we have $r \in \ZZZ(\KKK)$.

\smallskip

Assume first that $\ov{c} \neq c$. Let $a = (c+\ov{c}-2)(\ov{c}-c)^{-1}$, it commutes with $c$ and is such that $\ov{a}=-a$. Let $(x_n)_{n \in \NNN }$ be a sequence in $\ZZZ(\KKK)$ going to infinity. Let us define the marked apartment $f_n \in \Mor_{inj}({\cal A},{\cal I})$ by $f_n(\<e_i\>) = \<\widetilde{f_n}(e_i)\>$ for all $i \in \lb 1,4 \rb$, where
\beq &\widetilde{f_n}(e_1)=e_1, \hspace{0.5cm} \widetilde{f_n}(e_2) = e_2, \hspace{0.5cm} \widetilde{f_n}(e_3) = (e_1+e_2)(2^{-1}(1-a)\ov{c} -1)rx_n^2+e_3+e_5x_n& \\
& \mbox{ and } \widetilde{f_n}(e_4) = -(e_1c+e_2)2^{-1}(1+a)rx_n^2+e_4+e_5 x_n.&\eeq

Each of these vectors is isotropic : $q(\widetilde{f_n}(e_1))=0$, $q(\widetilde{f_n}(e_2))=0$,
\beq q(\widetilde{f_n}(e_3)) &=& (2^{-1}(1-a)\ov{c} -1)rx_n^2+(2^{-1}(1+a)c -1)rx_n^2+rx_n^2 \\
&=& 2^{-1}(\ov{c}+c-a(\ov{c}-c)-2)rx_n^2 = 0 \mbox{ and } \\
q(\widetilde{f_n}(e_4)) &=& -2^{-1}(1+a)rx_n^2-2^{-1}(1-a)rx_n^2+rx_n^2=0.\eeq

Furthermore $\varphi(\widetilde{f_n}(e_1),\widetilde{f_n}(e_2)) = \varphi(\widetilde{f_n}(e_2),\widetilde{f_n}(e_3)) = \varphi(\widetilde{f_n}(e_4),\widetilde{f_n}(e_1)) = 0$, and
$$ \varphi(\widetilde{f_n}(e_3),\widetilde{f_n}(e_4)) = (2^{-1}(1+a)c -1)rx_n^2-c2^{-1}(1+a)rx_n^2 + rx_n^2 = 0.$$
Hence $f_n$ is a marked apartment of ${\cal I}$. And it is clear that $(f_n)_{n \in \NNN}$ converges to $f$, so $f \in \ov{\Mor_{inj}({\cal A},{\cal I})}^g$.

If $\ov{c} = c$, then $f$ is a limit of plane-type quadripods whose cross-ratios do not belong to $\KKK^\sigma$, hence according to the previous case $f$ also belongs to $\ov{\Mor_{inj}({\cal A},{\cal I})}^g$.
\eit

\bigskip

It remains to show that the morphisms $t''_{12}$ and $t''_{23}$ belong to $\ov{\Mor_{inj}({\cal A},{\cal I})}^g$.

Firstly, the morphism $t''_{12}$ can be interpreted as a degenerate plane-type quadripod, and it is a limit of nondegenerate plane-type quadripods, so $t''_{12} \in \ov{\Mor_{inj}({\cal A},{\cal I})}^g$.

Secondly, the morphism $t''_{23}$ can be interpreted as a degenerate line-type quadripod, for which two opposite plane-type vertices coincide: it is symmetric. Then there exists a sequence $(f_n)_{n \in \NNN}$ of symmetric nondegenerate line-type quadripods which converges to $t''_{23}$. Hence each $f_n$ belongs to $\ov{\Mor_{inj}({\cal A},{\cal I})}^g$, and so does $t''_{23}$. 
\eit
\ep

This concludes the proof of Theorem~\ref{thmi:c2}.

\subsection{Rank $3$ : $\PGL(4)$}

Let $\KKK$ be a local field, and let ${\cal I}$ be the topological spherical building of complete flags of the vector space $V=\KKK^4$: it is the spherical building of $(\PGL_4,\KKK)$. Let ${\cal A}$ be the standard apartment in ${\cal I}$ which corresponds to the four canonical points $\{\< e_1\>,\<e_2\>,\<e_3\>,\<e_4\>\}$ of $\PPP(V)$. Denote by $\lb 1,4\rb_2$ the set of subsets of $\lb 1,4\rb$ of cardinality $2$.

If $f:{\cal A} \ra {\cal I}$ is a type-preserving simplicial morphism, we will denote
\beq \forall i \in \lb 1,4 \rb, & x_i = f(\<e_i\>) \\
\forall \{i,j\} \in \lb 1,4\rb_2, & \ell_{ij} = f(\<e_i,e_j\>) \\
\forall i \in \lb 1,4 \rb, & p_i = f(\<e_j \,:\, j \neq i\>),\eeq
where $\lb 1,4\rb_2$ denote the set of all subsets of $\lb 1,4 \rb$ with $2$ elements. Choosing a marked apartment of ${\cal I}$ is the same thing as choosing $4$ points, $6$ lines and $4$ planes in $\PPP(V)$,
$$ ((x_i)_{i \in \lb 1,4 \rb},(\ell_{ij})_{\{i,j\} \in \lb 1,4\rb_2},(p_i)_{i \in \lb 1,4 \rb}) \in \PPP(V)^4 \times \Gr_2(V)^6 \times \PPP(V^*)^4,$$
which satisfy the incidence conditions of a tetrahedron (see Figure~\ref{fig:tetraedre}),

\begin{multicols}{2}

\beq &\forall i,j \in \lb 1,4\rb \mbox{ distinct},&\\ &x_i \subset \ell_{ij}& \eeq
\vspace{-1cm}
\beq &\forall i,j,k \in \lb 1,4 \rb \mbox{ pairwise distinct},&\\ &\ell_{ij} \subset p_k& \eeq

\columnbreak

\begin{Figure}
\begin{center}
\includegraphics[height=3cm]{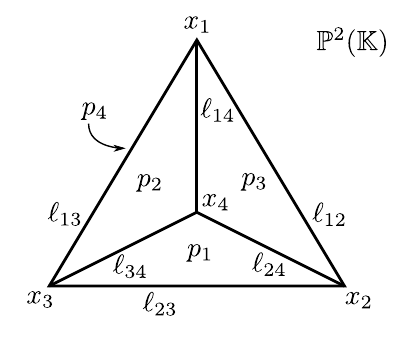}
\captionof{figure}{An apartment: a tetrahedron in $\PPP^2(\KKK)$}
\label{fig:tetraedre}
\end{center}
\end{Figure}
\end{multicols}


If $(f_n)_{n \in \NNN}$ is a sequence in $\Mor({\cal A},{\cal I})$, we will denote similarly
\beq \forall i \in \lb 1,4 \rb, & x^n_i = f_n(\<e_i\>) \\
\forall \{i,j\} \in \lb 1,4\rb_2, & \ell^n_{ij} = f_n(\<e_i,e_j\>) \\
\forall i \in \lb 1,4 \rb, & p^n_i = f_n(\<e_j \,:\, j \neq i\>).\eeq

If $f \in \Mor({\cal A},{\cal I})$ is such that there exists $\ell \in \Gr_2(V)$ such that $\forall \{i,j\} \in \lb 1,4\rb_2, \ell_{ij}=\ell$, we will say that $f$ is of type $(L)$ (see Figure~\ref{fig:typel}). In that case, we will say that $f$ is \df{symmetric} if $\Card \{x_1,x_2,x_3,x_4\} \leq 2$, or if $\Card \{p_1,p_2,p_3,p_4\} \leq 2$, or otherwise if there exists a projective isomorphism $s : \PPP(V) \simeq \PPP(V^*)$ such that for all $i \in \lb 1,4 \rb$, we have $s(x_i)=p_i$. This last condition is equivalent to asking that the cross-ratio of the four points $x_1,\ldots,x_4$ on the projective line $\ell$ is equal to the cross-ratio of the four planes $p_1,\ldots,p_4$ on the projective line $\ell^\perp$ of $\PPP(V^*)$.

Recall that if $x_1,x_2,x_3,x_4$ are points of the projective line $\KKK\PPP^1 \simeq \KKK \cup \{\infty\}$ such that $x_1$, $x_2$ and $x_3$ are pairwise distinct, the cross-ratio of $(x_1,x_2,x_3,x_4)$ is defined as $f(x_4) \in \KKK\PPP^1 \simeq \KKK \cup \{\infty\}$, where $f$ is the unique homography such that $f(x_1)=\infty$, $f(x_2)=0$ and $f(x_3)=1$. If $x_1,x_2,x_3,x_4$ are points of a projective line such that $\Card \{x_1,x_2,x_3,x_4\} \geq 3$, choose a double transposition $\sigma$ of $\lb 1,4 \rb$ such that $x_{\sigma(1)}$, $x_{\sigma(2)}$ and $x_{\sigma(3)}$ are pairwise distinct, then define the cross-ratio of $(x_1,x_2,x_3,x_4)$ to be the cross-ratio of $(x_{\sigma(1)},x_{\sigma(2)},x_{\sigma(3)},x_{\sigma(4)})$ (note that it is independent on the choice of $\sigma$). Hence the cross-ratio is a continous map from $\{(x_1,x_2,x_3,x_4) \in (\KKK\PPP^1)^4 \,:\, \Card \{x_1,x_2,x_3,x_4\} \geq 3\}$ to $\KKK \cup \{\infty\}$.

\begin{multicols}{2}
\hspace{1.5cm}
\begin{Figure}
\begin{center}
\includegraphics[height=2.5cm]{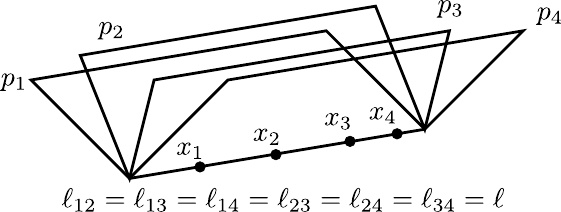}
\captionof{figure}{A morphism of type $(L)$}
\label{fig:typel}
\end{center}
\end{Figure}
\columnbreak
\begin{Figure}
\begin{center}
\includegraphics[height=4cm]{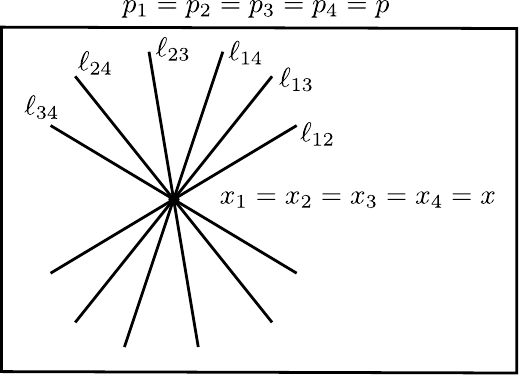}
\captionof{figure}{A morphism of type $(XP)$}
\label{fig:typexp}
\end{center}
\end{Figure}
\end{multicols}

If $f \in \Mor({\cal A},{\cal I})$ is such that there exists $x \in \PPP(V)$ and $p \in \PPP(V^*)$ such that for every $i \in \lb 1,4 \rb$, $x_i=x$ and $p_i=p$, we will say that $f$ is of type $(XP)$ (see Figure~\ref{fig:typexp}). In that case, we will say that $f$ is \df{symmetric} if there exists three pairwise intersecting $a,b,c \in \lb 1,4\rb_2$ such that $\ell_a=\ell_b=\ell_c$, or if there exists an involution $s$ of the projective line $\{\ell \in \Gr_2(V) \,:\, x \subset \ell \subset p\}$ such that for all $a \in \lb 1,4\rb_2$, we have $s(\ell_a)=\ell_{\lb 1,4\rb \bs a}$.

\bigskip

Notice that if $f$ is both of types $(L)$ and $(XP)$ (that is, if the image of $f$ is a complete flag), then $f$ is symmetric for both definitions.

\blem \label{lem:typeL} Let $f \in \Mor({\cal A},{\cal I})$ be of type $(L)$. Then $f$ belongs to $\ov{\Mor_{inj}({\cal A},{\cal I})}^g$ if and only if $f$ is symmetric. \elem

\bp Fix $f \in \ov{\Mor_{inj}({\cal A},{\cal I})}^g$ of type $(L)$, and let $(f_n)_{n \in \NNN}$ be a sequence of marked apartments which converges to $f$. Call $\ell$ the common line of $f$. Let us show that $f$ is symmetric. Fix a line $\ell_\infty$ in $\PPP(V)$ generic for $f$ and every $(f_n)_{n \in \NNN}$, that is, $\ell$ is such that $\ell \cap \ell_\infty = \emptyset$ and $\forall n \in \NNN, \forall \{i,j\} \in \lb 1,4\rb_2, \ell_{ij}^n \cap \ell_\infty = \emptyset$.

Fix $n \in \NNN$. Since $(x_i^n)_{i \in \lb 1,4 \rb}$ and $(p_i^n)_{i \in \lb 1,4 \rb}$ are projective frames of $\PPP(V)$ and $\PPP(V^*)$ respectively, the space of projective isomorphisms $s : \PPP(V) \simeq \PPP(V^*)$ such that for all $i \in \lb 1,4 \rb$ we have $s(x_i^n)=p_i^n$ has dimension $3$. The condition that $s(\ell_\infty)=\ell_\infty$ is given by $3$ independent linear homogeneous equations, hence there exists a unique isomorphism $s_n$ satisfying both properties.

Hence for all $n \in \NNN$, we have the equality between the two cross-ratios, the first one being in $\ell_\infty^\perp$ and the second one in $\ell_\infty$
$$[\<x^n_i,\ell_\infty\>]_{1 \leq i \leq 4}=[\<p^n_i \cap \ell_\infty\>]_{1 \leq i \leq 4}.$$
If $\Card \{x_1,x_2,x_3,x_4\} \leq 2$, or if $\Card \{p_1,p_2,p_3,p_4\} \leq 2$, then $f$ is symmetric. Otherwise, since the cross-ratio is continuous, taking the limit as $n$ goes on infinity gives
$$ [\<x_i,\ell_\infty\>]_{1 \leq i \leq 4}=[\<p_i \cap \ell_\infty\>]_{1 \leq i \leq 4}.$$
Since the maps $\ell_\infty^\perp \ra \ell : q \mapsto q \cap \ell$ and $\ell_\infty \ra \ell^\perp : y \mapsto \<y,\ell\>$ are isomorphisms of projective lines, we deduce that
$$[x_i]_{1 \leq i \leq 4}=[p_i]_{1 \leq i \leq 4},$$
hence $f$ is symmetric.


\bigskip

Conversely, fix $f \in \Mor({\cal A},{\cal I})$ of type $(L)$ and symmetric. Let us show that $f$ belongs to $\ov{\Mor_{inj}({\cal A},{\cal I})}^g$.
\ben

\item Assume first that $\Card \{p_1,p_2,p_3,p_4\} = 2$, for instance $p_1=p_2=p_3$. Assume first that $\Card \{x_1,x_2,x_3,x_4\} =4$. Up to precomposing $f$ by an element of the Weyl group and postcomposing $f$ by an element of $\PGL(V)$, we may assume that there exists $a \in \KKK \bs \{0,1\}$ such that
\beq &p_1 = p_2 = p_3 = \<e_1,e_2,e_3\>, p_4 = \<e_1,e_2,e_4\>,&\\
&x_1 = \<e_1\>, x_2 = \<e_2\>, x_3=\<e_1+e_2\> \mbox{ and } x_4=\<e_1a+e_2\>.& \eeq
Fix a sequence $(\alpha_n)_{n \in \NNN}$ in $\KKK \bs \{0\}$ converging to $0$. For all $n \in \NNN$, let $f_n$ be the marked apartment defined by
$$ x^n_1= \<e_1\>, x^n_2 = \<e_2\>, x^n_3 = \<e_1+e_2+e_4\alpha_n^2\> \mbox{ and } x^n_4 = \<e_1a+e_2+e_3\alpha_n\>.$$
The sequence $(f_n)_{n \in \NNN}$ converges to $f$. And if $\Card \{x_1,x_2,x_3,x_4\} < 4$, then $f$ is a limit of cases where the four points are distinct.

\item Assume $\Card \{x_1,x_2,x_3,x_4\} = 2$, then by duality we are in the previous case.

\item Assume that $\Card \{x_1,x_2,x_3,x_4\} = \Card \{p_1,p_2,p_3,p_4\} =1$. Then the image of $f$ is simply a complete flag in $V$, and it belongs to $\ov{\Mor_{inj}({\cal A},{\cal I})}^g$.

\item Assume that we are in none of the previous cases, then $\Card \{x_1,x_2,x_3,x_4\} \geq 3$ and $\Card \{p_1,p_2,p_3,p_4\} \geq 3$. Up to precomposing $f$ by an element of the Weyl group and postcomposing $f$ by an element of $\PSL(V)$, we may assume that there exists $a \in \KKK \cup \{\infty\}$ (the cross-ratio of $(x_1,x_2,x_3,x_4)$) such that
\beq &x_1 = \<e_1\>, x_2 = \<e_2\>, x_3=\<e_1+e_2\>, x_4=\<e_1a+e_2\>,&\\
&p_1 = \<e_1,e_2,e_3\>, p_2 = \<e_1,e_2,e_4\>, p_3 = \<e_1,e_2,e_3+e_4\> \mbox{ and } p_4 = \<e_1,e_2,e_3a+e_4\>.& \eeq
Fix a sequence $(\alpha_n)_{n \in \NNN}$ in $\KKK \bs \{0\}$ converging to $0$. For all $n \in \NNN$, let $f_n$ be the marked apartment defined by
$$ x^n_1= \<e_1\>, x^n_2 = \<e_2\>, x^n_3 = \<e_1+e_2+e_3a\alpha_n+e_4\alpha_n\> \mbox{ and } x^n_4 = \<e_1a+e_2+e_3a\alpha_n+e_4a\alpha_n\>.$$
Then the sequence $(f_n)_{n \in \NNN}$ converges to $f$. 
\een
\ep

\blem \label{lem:typeXP} Let $f \in \Mor({\cal A},{\cal I})$ be of type $(XP)$. Then $f$ belongs to $\ov{\Mor_{inj}({\cal A},{\cal I})}^g$ if and only if $f$ is symmetric. \elem

\bp Fix $f \in \ov{\Mor_{inj}({\cal A},{\cal I})}^g$ of type $(XP)$, and let $(f_n)_{n \in \NNN}$ be a sequence of marked apartments which converges to $f$. Call $x$ and $p$ the common point and plane of $f$. Let us show that $f$ is symmetric. Fix a point $x_\infty \in \PPP(V)$ generic for $f$ and every $(f_n)_{n \in \NNN}$, that is, $x_\infty$ is such that $x_\infty \not\in p$ and $\forall n \in \NNN, \forall i \in \lb 1,4\rb, x_\infty \not\in p_i^n$. Fix a line $\ell_\infty$ in $p$ that does not contain $x$.

Consider the continuous projection
\beq \pi : \PPP(V) \bs p & \lra & p \\
x' & \lmapsto & \<x',x_\infty\> \cap p. \eeq

Then the sequence $(g_n = \pi \circ f_n)_{n \in \NNN}$ converges to $\pi \circ f = f$.

Fix $n \in \NNN$ such that $\forall i \in \lb 1,4\rb, \pi(x_i^n) \not\in \ell_\infty$. Then according to Desargues'
involution theorem (see~\cite[Chap.~II, Theorem~31]{samuel_projective} or \cite[Chap.~1, Théorème~3.5.6]{perrin_projective}), there exists a (unique, projective) involution $s_n$ of $\ell_\infty$ such that for all $a \in \lb 1,4\rb_2$, we have $s_n(\pi(\ell_a^n) \cap \ell_\infty)=\pi(\ell_{\lb 1,4\rb \bs a}^n) \cap \ell_\infty$.

If there exists three pairwise intersecting $a,b,c \in \lb 1,4\rb_2$ such that $\ell_a=\ell_b=\ell_c$, then $f$ is symmetric. Assume this is not the case, then there exists a unique involution $s$ of $\ell_\infty$ such that $s(\ell_{12} \cap \ell_\infty)=\ell_{34} \cap \ell_\infty$ and $s(\ell_{13} \cap \ell_\infty) = \ell_{24} \cap \ell_\infty$. Since for all $a \in \lb 1,4\rb_2$, the sequence of points $(\pi(\ell_a^n))_{n \in \NNN}$ converges to $\ell_a$, we know that the sequence of involutions $(s_n)_{n \in \NNN}$ of $\ell_\infty$ converges to $s$, hence $s(\ell_{14} \cap \ell_\infty) = \ell_{23} \cap \ell_\infty$. So $f$ is symmetric.

\bigskip

Conversely, let $f \in \Mor({\cal A},{\cal I})$ be of type $(XP)$ and symmetric. Let us show that $f$ belongs to $\ov{\Mor_{inj}({\cal A},{\cal I})}^g$. Let us call $x$ and $p$ the common point and plane of $f$.

Assume first that there does not exist three pairwise intersecting $a,b,c \in \lb 1,4\rb_2$ such that $\ell_a=\ell_b=\ell_c$. Further assume that for all $a,b \in \lb 1,4 \rb_2$ such that $a \not\in \{b,\ov{b}\}$, we have $\{\ell_a,\ell_{\ov{a}}\} \cap \{\ell_b,\ell_{\ov{b}}\} \neq \emptyset$: up to precomposing $f$ by an element of the Weyl group, we may assume that $\ell_{12}=\ell_{13}$, $\ell_{34}=\ell_{23}$ and $\ell_{24}=\ell_{14}$. Since $f$ is symmetric, there exists a projective involution of the projective line $\{\ell \in \Gr_2(V) \,:\, x \subset \ell \subset p\}$ exchanging each pair of these three lines, so we have $\ell_{12}=\ell_{13}=\ell_{34}=\ell_{23}=\ell_{24}=\ell_{14}$: this contradicts the assumption.

So we can assume that there exist $a,b \in \lb 1,4 \rb_2$ such that $a \not\in \{b,\ov{b}\}$ and $\{\ell_a,\ell_{\ov{a}}\} \cap \{\ell_b,\ell_{\ov{b}}\} = \emptyset$. Up to precomposing $f$ by an element of the Weyl group, we may assume that $a=\{12\}$ and $b=\{13\}$. Assume further that $\ell_{12} \neq \ell_{34}$. Up to postcomposing $f$ by an element of $\PGL(V)$, we may assume that there exist $u,v,w \in \KKK \cup \{\infty\}$ such that $x=\<e_3\>$, $p=\<e_1,e_2,e_3\>$ and
\beq &\ell_{12} = \<e_1,e_3\>, \ell_{13} = \<e_1+e_2,e_3\>, \ell_{14} = \<e_1u+e_2,e_3\>,&\\
&\ell_{23} = \<e_1v+e_2,e_3\>, \ell_{24} = \<e_1w+e_2,e_3\> \mbox{ and } \ell_{34} = \<e_2,e_3\>.&\eeq
The assumption we made tells us that $w \not\in \{0,\infty\}$, and the symmetry condition tells us that there exists a projective involution $s$ of the projective line $\{\ell \in \Gr_2(V) \,:\, x \subset \ell \subset p\}$ such that for all $c \in \lb 1,4 \rb_2$ we have $s(\ell_c)=\ell_{\ov{c}}$. This involution is unique, it is defined by
\beq s : \{\ell \in \Gr_2(V) \,:\, x \subset \ell \subset p\} & \ra & \{\ell \in \Gr_2(V) \,:\, x \subset \ell \subset p\} \\
\ell=\<e_1z+e_2,e_3\> & \mapsto & \<e_1\frac{w}{z}+e_2,e_3\>.\eeq
Since $s(\ell_{14})=\ell_{23}$, we deduce that $\frac{w}{u}=v$.

\smallskip

Assume further that $u,v \not\in \{0,\infty\}$. Fix a sequence $(\alpha_n)_{n \in \NNN}$ in $\KKK \bs \{0\}$ going to $\infty$. For all $n \in \NNN$, let $f_n$ be the marked apartment defined by
\beq &x^n_1= \<e_3\>, x^n_2 = \<e_1(1-v)+e_2\alpha_n^{-1}+e_3\alpha_n\>,&\\
&x^n_3 = \<e_1+e_2+e_3\alpha_n\> \mbox{ and } x^n_4 = \<e_1+e_2u^{-1}+e_3\alpha_n+e_4\alpha_n^{-1}\>.&\eeq
Then if $a \in \{12,13,14\}$, it is easy to see that the sequence of lines $(\ell_a^n)_{n \in \NNN}$ converges to $\ell_a$. Every accumulation point of the sequence of lines $(\ell_{23}^n)_{n \in \NNN}$ contains the vector
$$ \liml_{n \ra +\infty} e_1v+e_2(1-\alpha_n^{-1}) = e_1v+e_2,$$
so it converges to $\ell_{23}$. And every accumulation point of the sequence of lines $(\ell_{24}^n)_{n \in \NNN}$ contains the vector
$$ \liml_{n \ra +\infty} e_1v+e_2(u^{-1}-\alpha_n^{-1})+e_4\alpha_n^{-1} = e_1v+e_2u^{-1} = (e_1w+e_2)u^{-1},$$
so it converges to $\ell_{24}$. And every accumulation point of the sequence of lines $(\ell_{34}^n)_{n \in \NNN}$ contains the vector
$$ \liml_{n \ra +\infty} e_2(1-u^{-1})-e_4\alpha_n^{-1} = e_2(1-u^{-1}),$$
so it converges to $\ell_{34}$.

Hence the sequence of marked apartments $(f_n)_{n \in \NNN}$ converges to $f$.

\smallskip

If for instance $u=0$ and $v=\infty$, then $f$ is a limit of previous cases with $u$ going to $0$ and $v=\frac{w}{u}$ going to infinity, so $f$ is a limit of marked apartments as well.

And if $\ell_{12} = \ell_{34}$, then $f$ is a limit of previous cases with $\ell_{12} \neq \ell_{34}$, so $f$ is a limit of marked apartments as well.

\bigskip

Assume now that there exist three pairwise intersecting $a,b,c \in \lb 1,4\rb_2$ such that $\ell_a=\ell_b=\ell_c$. Up to precomposing $f$ by an element of the Weyl group, we may assume that $\ell_{14}=\ell_{24}=\ell_{34}$. Fix $v \in \KKK \bs \{0\}$, and a sequence $(\alpha_n)_{n \in \NNN}$ in $\KKK \bs \{0\}$ going to $\infty$. For all $n \in \NNN$, let $f_n$ be the morphism of type $(XP)$ with common point $x=\<e_3\>$, common plane $p=\<e_1,e_2,e_3\>$ and the following lines:
\beq &\ell_{12} = \<e_1,e_3\>, \ell_{13} = \<e_1+e_2,e_3\>, \ell_{14} = \<e_1\alpha_n+e_2,e_3\>,&\\
&\ell_{23} = \<e_1v+e_2,e_3\>, \ell_{24} = \<e_1v\alpha_n+e_2,e_3\> \mbox{ and } \ell_{34} = \<e_2,e_3\>.&\eeq
According to the previous case, $f_n$ is a limit of marked apartments. And the sequence $(f_n)_{n \in \NNN}$ converges to the morphism of type $(XP)$ with the following lines:
\beq &\ell_{12} = \<e_1,e_3\>, \ell_{13} = \<e_1+e_2,e_3\>, \ell_{14} = \<e_2,e_3\>,&\\
&\ell_{23} = \<e_1v+e_2,e_3\>, \ell_{24} = \<e_2,e_3\> \mbox{ and } \ell_{34} = \<e_2,e_3\>.&\eeq
Hence we have shown that a generic morphism of type $(XP)$ such that $\ell_a=\ell_b=\ell_c$ is a limit of marked apartments, hence any morphism of type $(XP)$ with $\ell_a=\ell_b=\ell_c$ is also a limit of marked apartments.
\ep

\bthm \label{thm:sl4} Let ${\cal I}$ be the topological spherical building of $(\PGL_4,\KKK)$, and let ${\cal A}$ be an apartment of ${\cal I}$, then any morhism of $\Mor({\cal A},{\cal I})$ which is not of type $(XP)$ or $(L)$ belongs to $\ov{\Mor_{inj}({\cal A},{\cal I})}^g$. And a morphism of type $(XP)$ or $(L)$ belongs to $\ov{\Mor_{inj}({\cal A},{\cal I})}^g$ if and only if it is symmetric.
\ethm

\bp Fix $f \in \Mor({\cal A},{\cal I})$.
\bit
\item Assume that there exists $i \in \lb 1,4\rb$ such that $x_i \not\in p_i$, then $f$ is not of type $(L)$ nor $(XP)$. Let us assume that $i=1$. According to Theorem~\ref{thm:limites_a2} applied to $p_1$, there exists three sequences of points $(x^n_2)_{n \in \NNN}$, $(x^n_3)_{n \in \NNN}$ and $(x^n_4)_{n \in \NNN}$ of $p_1$ (which are not aligned for any $n \in \NNN$), which converge to $x_2$, $x_3$ and $x_4$ respectively, and such that the three sequence of lines $(\ell^n_{23})_{n \in \NNN}$, $(\ell^n_{34})_{n \in \NNN}$ and $(\ell^n_{42})_{n \in \NNN}$ they define converge to $\ell_{23}$, $\ell_{34}$ and $\ell_{42}$ respectively. For all $n \in \NNN$, let $f_n \in \Mor_{inj}({\cal A},{\cal I})$ be the marked apartment whose vertices are $x_1$, $x^n_2$, $x^n_3$ and $x^n_4$. Then the sequence $(f_n)_{n \in \NNN}$ converges to $f$.

\item Assume that there exists a partition $\lb 1,4 \rb = \{i,j\} \sqcup \{k,l\}$ such that $\ell_{ij} \cap \ell_{kl} = \emptyset$, then $f$ is not of type $(L)$ nor $(XP)$. Let us assume that $\{i,j\}=\{1,2\}$ and $\{k,l\}=\{3,4\}$. According to Proposition~\ref{pro:limites_a12} applied to $\ell_{12} \oplus \ell_{34}$, there exists four sequences of pairwise distinct points $(x^n_1)_{n \in \NNN}$, $(x^n_2)_{n \in \NNN}$, $(x^n_3)_{n \in \NNN}$ and $(x^n_4)_{n \in \NNN}$, which converge to $x_1$, $x_2$, $x_3$ and $x_4$ respectively, such that $x^n_1, x^n_2 \in \ell_{12}$ and $x^n_3, x^n_4 \in \ell_{34}$. For all $n \in \NNN$, let $f_n \in \Mor_{inj}({\cal A},{\cal I})$ be the marked apartment whose vertices are $x^n_1$, $x^n_2$, $x^n_3$ and $x^n_4$. Then the sequence $(f_n)_{n \in \NNN}$ converges to $f$.

\item Assume that there exists $p \in \PPP(V^*)$ such that for all $i \in \lb 1,4 \rb$ we have $p_i=p$, and that there is a point of $\{x_1,\ldots,x_4\}$ distinct from the three others. For instance $x_1 \not\in \{x_2,x_3,x_4\}$. According to Theorem~\ref{thm:limites_a2} applied to $p_1$, there exist three sequences of points $(x^n_2)_{n \in \NNN}$, $(x^n_3)_{n \in \NNN}$ and $(x^n_4)_{n \in \NNN}$ of $p_1$ (which are not aligned for any $n \in \NNN$), which converge to $x_2$, $x_3$ and $x_4$ respectively, and such that the three sequence of lines $(\ell^n_{23})_{n \in \NNN}$, $(\ell^n_{34})_{n \in \NNN}$ and $(\ell^n_{42})_{n \in \NNN}$ converge to $\ell_{23}$, $\ell_{34}$ and $\ell_{42}$ respectively. Choose a sequence $(x^n_1)_{n \in \NNN}$ in $\PPP(V) \bs p$ converging to $x_1$. For all $n \in \NNN$, let $f_n \in \Mor_{inj}({\cal A},{\cal I})$ be the marked apartment whose vertices are $x^n_1$, $x^n_2$, $x^n_3$ and $x^n_4$. Then the sequence $(f_n)_{n \in \NNN}$ converges to $f$.

\item Assume that there exists $p \in \PPP(V^*)$ such that for all $i \in \lb 1,4 \rb$ we have $p_i=p$, and that for instance $x_1=x_2 \neq x_3=x_4$. Choose a sequence $(x^n_1)_{n \in \NNN}$ in $\ell_{12} \bs \{x_1\}$ converging to $x_1$. For all $n \in \NNN$, let $f_n \in \Mor({\cal A},{\cal I})$ be the morphism which differs from $f$ only by $f_n(\<e_1\>)=x^n_1$. Then the sequence $(f_n)_{n \in \NNN}$ converges to $f$. According to the previous case, we know that $f_n \in \ov{\Mor_{inj}({\cal A},{\cal I})}^g$ for all $n \in \NNN$, so $f \in \ov{\Mor_{inj}({\cal A},{\cal I})}^g$.

\item Assume that there exists $p \in \PPP(V^*)$ such that for all $i \in \lb 1,4 \rb$ we have $p_i=p$, and that $x_1=x_2=x_3=x_4$. Then $f$ is of type $(XP)$, and this case is covered by Lemma~\ref{lem:typeXP}.

\item Assume that there exists $x \in \PPP(V)$ such that for all $i \in \lb 1,4 \rb$ we have $x_i=x$, then by duality between $\PPP(V)$ and $\PPP(V^*)$, we are in one of the three previous cases.

\begin{figure}[!h]
\begin{center}
\includegraphics[height=5cm]{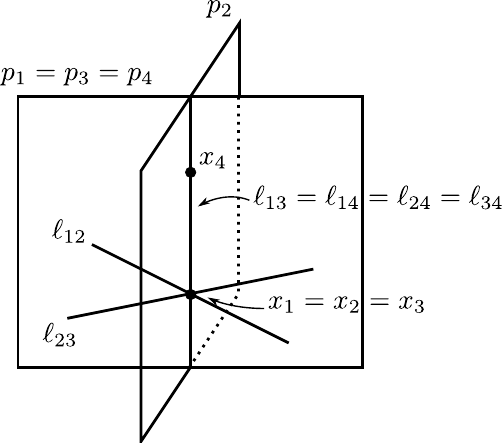}
\caption{The morphism in Lemma~\ref{lem:lastcase}}
\label{fig:lastcase}
\end{center}
\end{figure}

\item Assume that we are in none of the previous cases. Consider the intersection $\bigcap_{i \in \lb 1,4 \rb} p_i$ of the four planes : since $f$ is not a marked apartment, it is not empty. And since we are in none of the previous cases, the intersection is not a plane nor a point, so it is a line $\ell = \bigcap_{i \in \lb 1,4 \rb} p_i$. By duality between $\PPP(V)$ and $\PPP(V^*)$, we also know that the projective subspace spanned by the $(x_i)_{i \in \lb 1,4 \rb}$ is $\ell$. If for all $\{i,j\} \in \lb 1,4\rb_2$ we have $\ell_{ij}=\ell$, then $f$ is of type $(L)$, and this case is covered by Lemma~\ref{lem:typeL}. Otherwise assume that $\ell_{12} \neq \ell$, and let $x = \bigcap_{\{i,j\} \in \lb 1,4\rb_2} \ell_{ij}$: it is strictly included in $\bigcap_{i \in \lb 1,4 \rb} p_i = \ell$ and non-empty by the previous cases, so it is a point $x \in \PPP(V)$. Similarly, let $h = \<\ell_{ij} \,:\, \{i,j\} \in \lb 1,4\rb_2 \> \in \PPP(V^*)$.
\bit
\item If $x_1=x_2=x \not\in \{x_3,x_4\}$ for instance, choose a sequence $(x^n_2)_{n \in \NNN}$ in $\ell_{12} \bs \{x\}$ which converges to $x$. For all $n \in \NNN$, let $f_n \in \Mor({\cal A},{\cal I})$ be the morphism which differs from $f$ only by $f_n(\<e_2\>)=x^n_2$. Since $x^n_2 \not\in \ell$, according to the previous cases we know that $f_n \in \ov{\Mor_{inj}({\cal A},{\cal I})}^g$. As the sequence $(f_n)_{n \in \NNN}$ converges to $f$, we know that $f \in \ov{\Mor_{inj}({\cal A},{\cal I})}^g$.
\item If $x_1=x_2=x_3=x \neq x_4$ for instance, then by duality we may assume that $p_1=p_3=p_4=p \neq p_2$ as well, and up to precomposing $f$ by an element of the Weyl group and postcomposing $f$ by an element of $\PGL(V)$, this case is covered by the following Lemma~\ref{lem:lastcase} (see Figure~\ref{fig:lastcase}).

\eit
\eit
\ep

\blem \label{lem:lastcase} Let $f \in \Mor({\cal A},{\cal I})$ be defined by $x_1=x_2=x_3=\<e_1\>$, $x_4 = \<e_2\>$, $p_1=p_3=p_4=\<e_1,e_2,e_3\>$, $p_2 = \<e_1,e_2,e_4\>$, $\ell_{13}=\ell_{14}=\ell_{24}=\ell_{34}=\<e_1,e_2\>$, $\ell_{12}=\<e_1,e_3\>$ and $\ell_{23} = \<e_1,ae_2+be_3\>$, where $(a,b) \in \KKK^2 \bs \{(0,0)\}$. Then $f$ belongs to $\ov{\Mor_{inj}({\cal A},{\cal I})}^g$. \elem

\bp
Assume first that $a$ and $b$ are both non-zero. Fix a sequence $(\alpha_n)_{n \in \NNN}$ in $\KKK \bs \{0\}$ going to $0$. For all $n \in \NNN$, let $f_n$ be the marked apartment defined by
$$ x^n_1= \<e_1\>, x^n_2 = \<e_1+e_3b\alpha_n\>, x^n_3 = \<e_1-e_2a\alpha_n+e_4\alpha_n^2\> \mbox{ and } x^n_4 = \<e_2\>.$$
Then it is easy to see that for all $i \in \lb 1,4 \rb$, the sequence of points $(x_i^n)_{n \in \NNN}$ converges to $x_i$ and the sequence of planes $(p_i^n)_{n \in \NNN}$ converges to $p_i$. And if $a \in \{12,13,14,24,34\}$, it is easy to see that the sequence of lines $(\ell_a^n)_{n \in \NNN}$ converges to $\ell_a$. Every accumulation point of the sequence of lines $(\ell_{23}^n)_{n \in \NNN}$ contains the vector
$$ \liml_{n \ra +\infty} e_2a\alpha_n +e_3b\alpha_n - e_4\alpha_n^2= e_2a+e_3b,$$
so it converges to $\ell_{23}$. Hence the sequence of marked apartments $(f_n)_{n \in \NNN}$ converges to $f$.

If $a=0$ or $b=0$, then $f$ is a limit of similar morphism for which we have $a$ and $b$ both non-zero, so $f$ belongs to $\ov{\Mor_{inj}({\cal A},{\cal I})}^g$.
\ep


\begin{thebibliography}{GKVMW12}

\bibitem[AB08]{abramenko_brown}
Peter Abramenko and Kenneth~S. Brown.
\newblock {\em {Buildings. Theory and Applications.\!}}, volume 248 of {\em
  Grad. Texts in Maths}.
\newblock Springer, 2008.

\bibitem[Bae52]{baer}
Reinhold Baer.
\newblock {\em {Linear algebra and projective geometry.}}
\newblock {Pure and Applied Mathematics, 2., New York: Academic Press Inc.},
  1952.

\bibitem[BGS85]{bgs}
Werner Ballmann, Mikhael Gromov, and Viktor Schroeder.
\newblock {\em {Manifolds of nonpositive curvature\!}}
\newblock {Progr.~Math. {\bf 61}}. {Birkh\"auser}, 1985.

\bibitem[BH99]{bridson_haefliger}
Martin~R. Bridson and Andr\'e Haefliger.
\newblock {\em {Metric \,spaces \,of \,non-positive \,curva\-ture}}, volume
  {319\!} of {\em {Grund.~math.~Wiss.}}
\newblock {Springer}, 1999.

\bibitem[Bou59a]{bourbaki_algebre9}
Nicolas Bourbaki.
\newblock {\em {Algèbre. Chapitre IX}}.
\newblock {Hermann}, {1959}.

\bibitem[Bou59b]{bourbaki_integration8}
Nicolas Bourbaki.
\newblock {\em {Int\'egration. Chapitre VIII}}.
\newblock {Hermann}, 1959.

\bibitem[Bou97]{bourdon}
Marc Bourdon.
\newblock {Immeubles hyperboliques, dimension conforme et rigidit\'e de
  Mostow}.
\newblock {\em {Geom. Funct. Anal.}}, 7:245--268, 1997.

\bibitem[BS87]{burns_spatzier}
Keith Burns and Ralf Spatzier.
\newblock {On topological Tits buildings and their classification}.
\newblock {\em {Publ. Math. Inst. Hautes \'Etudes Sci.}}, 65:5--34, 1987.

\bibitem[CDP07]{cdp}
Gilles Courtois, Françoise Dal'bo, and Fr\'ed\'eric Paulin.
\newblock {\em {Sur la dynamique des groupes de matrices et applications
  arithm\'etiques\!}}
\newblock {Journ\'ees math\'ematiques X-UPS 2007}. {Les \'{E}ditions de
  l'\'{E}cole Polytechnique}, 2007.

\bibitem[CEG87]{ceg}
R.D. Canary, D.B.A Epstein, and P.L. Green.
\newblock {\em {Notes on notes of Thurston}}.
\newblock in "Analytical and geometric aspects of hyperbolic space",
  Lond.~Math.~Soc., Lect.~Notes~Series {\bf 328}, 1987.

\bibitem[Cha50]{chabauty}
Claude Chabauty.
\newblock {Limite d'ensembles et g\'eom\'etrie des nombres}.
\newblock {\em {Bull. Soc. Math. France}}, 78:143--151, 1950.

\bibitem[dlH08]{harpe_chabauty}
Pierre de~la Harpe.
\newblock {Spaces of closed subgroups of locally compact groups}.
\newblock {arXiv:0807.2030v2}, 2008.

\bibitem[FM94]{fulton_macpherson}
William Fulton and Robert MacPherson.
\newblock {A compactification of configuration spaces}.
\newblock {\em {Annals of Math.}}, 139:183--225, 1994.

\bibitem[GJT98]{guivarch}
Yves Guivarc'h, Lizhen Ji, and J.~C. Taylor.
\newblock {\em {Compactifications of symmetric spaces\!}}
\newblock {Progr.~Math. {\bf 156}}. {Birkh\"auser}, 1998.

\bibitem[GKK95]{kramer_flag_homogeneous}
T.~Grundh{\"o}fer, N.~Knarr, and L.~Kramer.
\newblock Flag-homogeneous compact connected polygons.
\newblock {\em Geom. Dedicata}, 55(1):95--114, 1995.

\bibitem[GKK00]{kramer_flag_homogeneous_II}
T.~Grundh{\"o}fer, N.~Knarr, and L.~Kramer.
\newblock Flag-homogeneous compact connected polygons. {II}.
\newblock {\em Geom. Dedicata}, 83(1-3):1--29, 2000.
\newblock Special issue dedicated to Helmut R. Salzmann on the occasion of his
  70th birthday.

\bibitem[GKVMW12]{kramer_compact_moufang_buildings}
Theo Grundh{\"o}fer, Linus Kramer, Hendrik Van~Maldeghem, and Richard~M. Weiss.
\newblock Compact totally disconnected {M}oufang buildings.
\newblock {\em Tohoku Math. J. (2)}, 64(3):333--360, 2012.

\bibitem[GP01]{gaboriau_paulin}
Damien Gaboriau and Fr\'ed\'eric Paulin.
\newblock {Sur les immeubles hyperboliques}.
\newblock {\em Geom. Dedicata}, 88:153--197, 2001.

\bibitem[Hae10a]{haettel_m2}
Thomas Haettel.
\newblock {Compactification de Chabauty des espaces sym\'etriques de type non
  compact}.
\newblock {\em {J. Lie Theory}}, 20:437--468, 2010.

\bibitem[Hae10b]{RxZ}
Thomas Haettel.
\newblock {L'espace des sous-groupes ferm\'es de $\mathbb{R} \times
  \mathbb{Z}$}.
\newblock {\em {Algebr. Geom. Topol.}}, 10:1395--1415, 2010.

\bibitem[Hae13]{haettel_chabauty_cartan}
Thomas Haettel.
\newblock {Compactification de Chabauty de l'espace des sous-groupes de Cartan
  de $\SL_n(\mathbb{R})$}.
\newblock {\em {Math. Z.}}, 274:573--601, 2013.

\bibitem[Har77]{harvey_discrete}
W.~J. Harvey.
\newblock {\em {Spaces of discrete groups}}.
\newblock in "Discrete groups and automorphic functions". 1977.

\bibitem[IM05a]{iliev_manivel_severi}
Atanas Iliev and Laurent Manivel.
\newblock {Severi varieties and their varieties of reductions}.
\newblock {\em {J. Reine Angew. Math.}}, 585:93--139, 2005.

\bibitem[IM05b]{iliev_manivel}
Atanas Iliev and Laurent Manivel.
\newblock {Varieties of reduction for ${\mathfrak{gl}_n}$}.
\newblock In {C. Ciliberto, A.V. Geramita, B. Harbourne, R. Miro-Roig, K.
  Ranestad}, editor, {\em {Projective Varieties with Unexpected Properties}},
  pages 287--316. Walter de Gruyter, 2005.

\bibitem[Ji06]{ji_buildings}
Lizhen Ji.
\newblock {Buildings and their applications in geometry and topology}.
\newblock {\em {Asian J. Math.}}, 10:11--80, 2006.

\bibitem[Kas12]{kassel_proper_actions}
Fanny Kassel.
\newblock {Deformation of proper actions on reductive homogeneous spaces}.
\newblock {\em {Math. Ann.}}, 353:599--632, 2012.

\bibitem[Kob09]{kobayashi_exposition}
Toshiyuki Kobayashi.
\newblock {On discontinuous group actions on non-Riemannian homogeneous
  spaces}.
\newblock {\em {Sugaku Expositions}}, 22:1--19, 2009.

\bibitem[LBG11a]{lebarbier_ex}
Michaà«l Le~Barbier~Grünewald.
\newblock {Examples of varieties of reductions of small rank}.
\newblock {\texttt{http://www.uni-bonn.de/$\sim$mlbg/public/michi-redex.pdf}},
  {2011}.

\bibitem[LBG11b]{lebarbier}
Michaà«l Le~Barbier~Grünewald.
\newblock {The variety of reductions for a reductive symmetric pair}.
\newblock {\em {Transform. Groups}}, 16:1--26, 2011.

\bibitem[OW02]{OW}
H.~Oh and D.~Witte.
\newblock {Compact Clifford-Klein forms of homogeneous spaces of $\SO(2,n)$}.
\newblock {\em {Geom. Dedicata}}, 89:25--57, 2002.

\bibitem[Per12]{perrin_projective}
Daniel Perrin.
\newblock {G\'eom\'etrie projective plane et applications aux g\'eom\'etries
  euclidienne et non euclidiennes}.
\newblock {http://www.math.u-psud.fr/\textasciitilde
  perrin/Livre\_de\_geometrie\_projective.html}, 2012.

\bibitem[Sam88]{samuel_projective}
Pierre Samuel.
\newblock {\em {Projective geometry}}.
\newblock {Undergrad. Texts in Math.} {Springer}, 1988.

\bibitem[Ser68]{serre_corpslocaux}
Jean-Pierre Serre.
\newblock {\em Corps locaux}.
\newblock Publ. Univ. Nancago, 1968.

\bibitem[Tit74]{tits_spherical_type}
Jacques Tits.
\newblock {\em {Buildings of spherical type and finite BN-pairs\!}}
\newblock {Lect. Notes in Math. {\bf 386}}. {Springer}, 1974.

\bibitem[TW02]{moufang}
Jacques Tits and Richard~M. Weiss.
\newblock {\em Moufang polygons}.
\newblock Springer Monogr. Math. Springer-Verlag, Berlin, 2002.

\end{thebibliography}
\end{document}